\documentclass[12pt]{article}
\newcommand{\assign}{:=}

\newcommand{\comment}[1]{}
\newcommand{\complexes}{{\bf C}}

\renewcommand{\div}{{\mathop{\rm div}\nolimits}}
\newcommand{\dist}{\mathop{\rm dist}\nolimits}

\newcommand{\grad}{\nabla}

\renewcommand{\Im}{\mathop{\rm Im}\nolimits}

\newcommand{\naturals}{{\bf N}}

\newcommand{\qed}{\vrule height 6pt depth 0pt width 3 pt}

\newcommand{\reals}{{\bf R}}
\renewcommand{\Re}{\mathop{\rm Re}\nolimits}

\newcommand{\supp}{\mathop{\rm supp}\nolimits }

\newcommand{\BibTeX}{{\rm B\kern-.05em{\sc i\kern-.025em b}\kern-.08em     
    T\kern-.1667em\lower.7ex\hbox{E}\kern-.125emX}}

\newcommand{\average}{-\!\!\!\!\!\!\int}
\newcommand{\e}{\epsilon}
\newcommand{\p}{\partial}
\newcommand{\Om}{\Omega}

\newcommand{\ball}[2]{B_{#2}(#1)}
\newcommand{\sball}[2]{\Delta_{#2}(#1)}
\newcommand{\nontan}[1]{#1^*}

\renewcommand{\grad}{\nabla}

\newcommand{\sobolev}[2]{W^{#1,#2}}

\newenvironment{proof}[1][Proof]{\begin{trivlist}\item[\hskip \labelsep
{\it #1. }]}{\hfill \qed \goodbreak \end{trivlist}}
\newenvironment{remark}[1][Remark]{\begin{trivlist}\item[\hskip \labelsep
{\em  #1. }]}{\goodbreak \end{trivlist}}

\newtheorem{theorem}{Theorem}[section]
\newtheorem{proposition}[theorem]{Proposition}

\newtheorem{lemma}[theorem]{Lemma}

\begin{document}

\title{The mixed problem in $L^p$ for some two-dimensional Lipschitz
domains}

\author{
Loredana Lanzani$^*$ \\ Department of Mathematics \\ University of
Arkansas \\ Fayetteville, Arkansas \\ lanzani@caccioppoli.uark.edu   \and 
Luca Capogna$^*$ \\ Department of Mathematics \\ University of Arkansas \\
Fayetteville, Arkansas \\
capogna@lagrange.uark.edu
\and 
Russell M. Brown
\footnote{Research supported, in part, by the National
Science Foundation.} 
\\ Department of Mathematics 
\\ University of Kentucky
\\ Lexington, Kentucky 
\\ russell.brown@uky.edu
}

\date{}
\maketitle

{
\abstract{  We consider the mixed problem, 
$$\left\{
\begin{array}{ll}
\Delta u = 0 \qquad & \mbox{in }  \Omega \\
\frac{\partial u }{\partial \nu} = f_N  \qquad  & \mbox{on } $N$ \\
u =f_D  \quad &  \mbox{on } D
\end{array}
\right. 
$$
in a class of Lipschitz graph domains in two dimensions with Lipschitz
constant at most 1. 
We suppose the Dirichlet data, $f_D$ has one derivative in $L^p(D)$ of
the boundary and the Neumann data is 
in $L^p(N)$. We find a $p_0>1$ so  that for $p$ in an interval
$(1,p_0)$, we may find  a unique solution to the mixed problem and the
gradient of the  solution lies in $L^p$. 
}

}

\section{Introduction}
\setcounter{equation}{0}
The goal of this paper is to study the mixed problem (or Zaremba's
problem)  for Laplace's
equation in certain two-dimensional domains when the Neumann data
comes from $L^p$ and the Dirichlet problem has one derivative in
$L^p$. We will consider both $L^p$ with respect to arc-length, $ d\sigma$ and
also $L^p(w\,d \sigma)$ where the weight $w$  is of the form $ w(x)=
|x|^\epsilon$. 
By the {\em{mixed problem for}} $L^p(w\, d\sigma)$, 
we mean the following
boundary value problem 
\begin{equation}\label{Mixed} \left\{
\begin{array}{ll}
\Delta u = 0 \qquad & \mbox{in } \Omega \\
u= f_D \qquad & \mbox{on } D \\
\frac { \partial u }{ \partial \nu } =
 f_N 
& \mbox{on } N \\
(\nabla u )^* \in   L^p(w\, d\sigma ) 
\end{array}
\right. 
\end{equation}
Here, $(\nabla u )^*$ is the non-tangential maximal function of the
gradient (see the definition in (\ref{defNT})). 
The domain $ \Omega$ will be a Lipschitz graph domain. Thus, $ \Omega
= \{ (x_1, x_2) : x_2> \phi(x_1)\}$ where $ \phi : \reals \rightarrow
\reals $ is a Lipschitz function with $\phi(0) =0 $.  We will call
such domains       {\em standard graph domains. }  
The sets $D$ and $N$ satisfy $ D\cup N = \partial \Omega$ and $ D\cap
N = \emptyset$. Also, we want $D$ to be open, so that it supports
Sobolev spaces. 
Here, we will generally assume that $f_D$ is a function with one
derivative in $L^p(D, w\,  d\sigma)$ and that $f_N$ is  in $L^p( N, w \,
d\sigma )$. 
Our goal is to find conditions on the data, the exponent $p$ and the
weight  $w$ so that the problem
(\ref{Mixed}) has a unique solution. This paper continues the work of
Sykes and Brown  \cite{RB:1994b,JS:1999,SB:2001} to provide a partial
answer to problem 3.2.15 from Kenig's CBMS lecture notes
\cite{CK:1994}. Our goal is to obtain $L^p$-results in a class of
domains that is not included in the domains studied by Sykes
\cite{JS:1999}.

We recall  a few works which study the regularity of the 
mixed problem. In \cite{AK:1982}, Azzam and Kreyszig established that the mixed
problem has a solution in $C^{2+\alpha}$ in bounded two-dimensional
domains provided that $N$ and $D$ meet at a sufficiently small angle.
Lieberman \cite{GL:1986} also gives conditions  which imply that the
solution is H\"older continuous.
In addition, much effort has been devoted to problems in polygonal
domains, see the monograph of Grisvard \cite{PG:1985}.
% and Mazya and Plamenevskii \cite{MP:????} among others. 
Savar\'e
\cite{GS:1997} finds solutions in the 
Besov space $B^{ 3/2}_{2,\infty}$  in smooth domains.  This positive
result fits quite nicely with the example below which shows that there
is a 
solution whose gradient just misses having non-tangential maximal function in
$L^2 ( d\sigma)$. It is known that solutions for which the
non-tangential maximal function of the gradient is in $L^2(d\sigma)$
also belong to the Besov space $B^{3/2}_{2,2}$ (see the article of
Fabes \cite{EF:1988}).

If we recall the standard tools for studying boundary value problems
in Lipschitz domains, we see that the mixed problem presents an
interesting technical challenge. 
On the one hand, the starting point for many results on boundary
value problems in Lipschitz domains is the Rellich
identity (see Jerison and Kenig \cite{JK:1982c}, for example). This
remarkable identity 
provides estimates at the boundary for derivatives of a 
harmonic function  in  $L^2$. On the other
hand, simple and well-known 
examples show that the mixed problem is not solvable in  
$L^2(d\sigma)$ for  smooth domains. Let us recall an example in the upper
half-space, 
$\reals^2 _{+} = \{(x_1,x_2) : x_2 >0\}$.  We let  
$N=\{(x_1,0), \ x_1>0\}$, $D=\{(x_1,0), \ x_1<0\}$, and we consider
the harmonic function
$$u(z)=Re\ (\sqrt{z}-\sqrt{z+i}), \qquad z=x_1+i x_2.$$
It is easy to see that we have 
\begin{eqnarray*}
\bigg|\frac{\p u}{\p \nu}\bigg| &\le& C(1+|z|)^{-3/2} \quad \mbox{on } \ N,\\
\bigg|\frac{d u}{d\sigma}\bigg| &\le& C(1+|z|)^{-3/2} \quad \mbox{on }\ D,
\end{eqnarray*}
but 
$$|\nabla u|\ge \frac{C}{\sqrt{|z|}} \qquad \mathrm{ for } \quad |z|\le 1/2.$$
 Hence, we do not have $(\nabla u)^*\in L^2_{loc}(\reals,d\sigma)$.
On the other hand we have $(\nabla u)^*\in L^p_{loc}(\reals,d\sigma)$ for 
all $1\le p<2$.

We detour around this problem by establishing weighted estimates in
$L^2$ using the Rellich identity. This relies on an observation of
Luis Escauriaza \cite{LE:2000} that a Rellich identity holds when the  components
 of the vector field
are, respectively, the real and imaginary part of a holomorphic function.
Then, we imitate the arguments of Dahlberg and Kenig \cite{DK:1987}
to establish
Hardy-space estimates (with a different weight). The weights are
chosen so that interpolation will give us unweighted $L^p $ as an
intermediate space.   
%We will occasionally use the machinery of $A_p$-weights. 
The weights we consider will  be of the form 
$|x|^\epsilon$ restricted to  the boundary of a Lipschitz graph domain. 
Earlier work of Shen \cite{ZS:2005} gives a different approach to the study of
weighted estimates for the Neumann and regularity  problems when the
weight is a power. 

In the result below and
throughout this paper, we assume that   $\Omega$ is a standard Lipschitz
graph
domain and that 
\begin{equation}\label{EDN}
 D= \{ ( x_1, \phi(x_1) ): x_1 < 0 \}\, ,\qquad
 N = \{ ( x_1, \phi(x_1) ) : x_1
\geq 0\}.
\end{equation}
 We will call $ \Omega$ with $N$ and $D$ as defined above, a
{\em standard Lipschitz graph domain for the mixed problem}.
 The {\em Lipschitz constant} of the domain is defined to be 
the quantity:
\begin{equation}\label{Lipconst}
M= \|\phi '\|_{L^{\infty}(\reals)}
\end{equation}  
Our main result is the following. 
\begin{theorem} \label{maintheorem}  Let $ \Omega$ be a standard Lipschitz graph domain 
 for the mixed problem, with
  Lipschitz constant $M$ less than 1. There exists $p_0= p_0(M)$
 so that for $1<p<p_0$,  if
 $f_N  \in L^p(N,d\sigma)$ and $ df_D/d\sigma \in L^p(D, d\sigma)$,
 then  
  the mixed problem for 
  $L^p(d\sigma)$ has a unique solution. The solution satisfies
\begin{equation}\label{EMPLP}
\| ( \nabla u )^ *\|_{ L^ p ( d\sigma)}\leq C(p,M)
 \left( \| f_N
  \|_{L^p(N,d\sigma)} +
 \left\| \frac{d f_D}{d\sigma}\right\| _{L^ p (D,d\sigma)}\right).
\end{equation} 
\end{theorem}

\section{Preliminaries}
\setcounter{equation}{0}

In this section, we prove uniqueness for the weighted $L^p$-Neumann and
regularity problems in a Lipschitz graph domain, see (\ref{Neumann}) and
(\ref{Regular}). The proofs are based on Rellich-type estimates (Proposition
\ref{RellichIneq}).
Here and in the sequel, we let $\nu$ denote the outer unit
 normal to $\partial \Omega$.
\comment{
\begin[L. Escauriaza]{lemma}
%\label{Mixed} 
\label{LuisE}  If $ \Delta u = 0$ in $\Omega$ and $
  \alpha _1 + i\alpha _2$ is 
holomorphic in $ \Omega$, then with $\alpha = ( \alpha _1 , \alpha
_2)$, we have 
$$
\div  (|\nabla u |^2 \alpha - (2 \alpha \cdot \grad u) \grad  u ) = 0.
$$
\end{lemma}
\comment{Works in C^n}
\begin{proof} The proof is a calculation. We omit the
details. \end{proof}

Using this Lemma and standard manipulations as in \cite{JK:1982c}, we
will obtain {\em a priori } estimates for our boundary value
problems. 

}
We recall that a {\em bounded} Lipschitz domain is a bounded
 domain whose boundary is parameterized by (finitely many) Lipschitz
 graphs. 
We begin our development with an observation that
we learned from Luis 
Escauriaza \cite{LE:2000}. 
\begin{lemma}[L. Escauriaza] \label{LuisE} Suppose $ \Omega $ is a bounded Lipschitz
domain and $ u$ and 
$\alpha =  (\alpha _1 , \alpha _2)$ are smooth in a neighborhood of $ \bar
\Omega$ with $ \Delta u =0$ and $ \alpha_1 + i \alpha _2 $
holomorphic.
Then, we have 
$$
\int_ {\partial \Omega } | \nabla u |^2 \alpha \cdot \nu - 2 
\frac {\partial u } {\partial \alpha} \frac { \partial  u }
{\partial \nu }   \, d\sigma = 0.
$$
\end{lemma}
\begin{proof}
A calculation shows that 
$\div  (|\nabla u |^2 \alpha - (2 \alpha \cdot \grad u) \grad  u ) = 0$.
Thus, the lemma is an immediate consequence of the divergence
theorem. 
\end{proof}

Next, we recall Carleson measures and the fundamental property of
these measures. The applications  we have in mind are  simple, but the
use of Carleson measures will allow us to appeal to a well-known
geometrical argument rather than invent our own.

Let $ \sigma$ be a  measure on the boundary of $ \Omega$.
A measure $\mu$ on $ \Omega $ is said to be a  { \em Carleson
measure with respect to $\sigma$ } if there is a constant $A$ so that
for  each $ x\in \partial
\Omega$ and each  $r>0$, we  have 
$$
\mu (\ball x r  \cap  \Omega ) \leq A  \sigma( \sball x r ).
$$
Here, we are using $ \ball x r$ to denote the ball (or disc)  in
$\reals ^2$ with center $x$ and radius $r$ and  we also use 
$\sball x r = \ball x r \cap \partial \Omega$ to denote a ball on the
boundary of $ \Omega$.   
%%Delete ball? 

Before we can state the next result,  we need a few
definitions.  
For $ \theta_0>0$,
we define $ \Gamma(0) = \{ re^{ 
i\theta} : r> 0, \ | \theta - \pi/2| < \theta _0\}$ to be the sector
at $0$ with vertical axis and 
 opening $ 2 \theta _0$. Then, we put $\Gamma (x) = x+
\Gamma(0)$, $x\in \partial \Omega$. 
If $ \Omega $ is a Lipschitz graph domain with constant
$M$, then we have that $\Gamma(x)$ defines a non-tangential approach
region provided $ \theta _0 < \pi/2 - \tan ^{ -1} (M)$.  We fix  such
a 
$ \theta _0$ and if $v$ is a function defined on $ \Omega$,
we  define the {\em non-tangential maximal function of} $v$, $\nontan{v} $ on
$\partial \Omega$  by 
\begin{equation}\label{defNT}
\nontan{v} (x) = \sup _{ y \in \Gamma (x) } |v (y)|\,
 ,\quad x\in \partial \Omega\, .
\end{equation}
We also use these sectors to define restrictions to the boundary in the
sense of non-tangential limits. For a function $v$ in $ \Omega$, we
define the restriction of $v$ to the boundary by 
\begin{equation}\label{Enontlim}
v(x) =\lim _{\Gamma(x) \ni  y\rightarrow x } v(y)\, , 
\quad x\in \partial \Omega
\end{equation}
provided the limit exists. 
Finally, we recall that a measure $ \sigma $ is a {\em doubling 
measure } if there is a constant $C$ so that for all $r>0$ we have $
\sigma( \Delta_{2r}(x) ) \leq C \sigma ( \Delta _r (x))$.   

\begin{proposition}  \label{Carleson} If $\tau $ is a doubling
measure on $ \partial \Omega$  and $\mu $ is a Carleson measure with
respect to $ \tau $ 
with constant $A$, then there exists a constant $C>0$ so that 
$$
\left|\int_{ \Omega} v \, d\mu \right| \leq CA  \int_{\partial \Omega}
\nontan{v} \, d\tau . 
$$
\end{proposition}

This  result is well-known. The proof in Stein
\cite[pp.~58--60]{ES:1993} easily generalizes from Lebesgue measure to
doubling measures. 

A simple example of a Carleson measure that will be useful to us is
the following.   With $ R>0$ and $ \epsilon > -1$, define  $ d \sigma
_ {\epsilon}$ and  $d \mu_\epsilon  $ by 
\begin{equation}\label{Ecarleson}
 d \sigma _ {\epsilon} (x)  = |x|^\epsilon d\sigma(x)\, ,\
 x\in \partial \Omega\quad  ;  \quad
d \mu_\epsilon (y) = \frac {|y|^\epsilon}  { R} \chi _{ \ball  0 R } (y) dy\, ,
\  y\in \Omega, 
\end{equation} 
where we are using  $d\sigma$ for arc-length
 measure on $\partial \Omega$ and $dy$ for area measure in the plane. 
It is not hard to see that $ \sigma _ \epsilon $ is a  doubling
measure on $\partial \Omega$  (see Lemma \ref{measure-property}) and
$\mu_\epsilon $ is a 
Carleson measure with respect to $\sigma_\epsilon $.

\begin{lemma} [Rellich Identity]  \label{RellichIdentity}
Given  $\epsilon > -1$ and $a \in \complexes$, we define 
$ \alpha (z) =  ( \Re (a z^ \epsilon) , \Im ( a z^ \epsilon ))$. 
Here, $z=x_1+ix_2$. Let $\sigma_\epsilon$ be as in
(\ref{Ecarleson}). If $u$ is harmonic in $ \Omega$ and  $\nontan {(\grad u) }
\in L^2 ( \sigma_\epsilon)$, then we have 
$$
\int _{ \partial \Omega } | \grad u | ^2 \alpha \cdot \nu - 2 \alpha
\cdot \nabla u \frac { \partial u }{ \partial \nu } \, d\sigma =0. 
$$
\end{lemma}
\begin{remark} Here and in the sequel, we let  $\nabla u (x)$
 denote the non-tangential limit of
$\nabla u$ at $x \in \partial\Omega$, see (\ref{Enontlim}). 
It is well known that, with the assumptions of
Lemma \ref{RellichIdentity}, the non-tangential limit
of $\nabla u$ exists a.e.
 $x \in \partial\Omega$, 
%since, with the assumptions of 
%Lemma \ref{RellichIdentity}, we have that the complex derivative of $u$,
%$\partial u$ (see (\ref{ECderiv})) is in the Hardy space
% $H^2(\partial \Omega)$,
see Dahlberg \cite{BD:1977} or Jerison and Kenig \cite{JK:1982d}.  
\end{remark} 
\begin{proof} We introduce a cut-off function $ \eta _R$ where $
  \eta_R (y)=1 $ if $ y\in B_R(0) $,
 $ \eta _R(y) = 0$ if $|y| >2R$, and $ | \grad \eta_R | \leq C/R$.
 For $\tau >0$, we define a  translate of $u$, $u_\tau$, by 
$u_ \tau(y) = u (y + \tau e_2)$, where $e_2 = (0,1)$.
Note that $u_\tau$ is smooth in a neighborhood of $\overline{\Omega}$.
 Since $\Omega$ is a graph domain and $0 \in \partial\Omega$
it follows that  $az^{\epsilon}$ is
holomorphic in $\Omega$.  Thus, we may apply Lemma \ref{LuisE}  to 
$u_\tau$ and $\alpha$.  The divergence theorem now  yields  
\begin{eqnarray*}
\int_{\partial \Omega} \left(  | \nabla u_\tau |^2 \alpha \cdot \nu - 2 \alpha
\cdot \nabla u _ \tau \frac { \partial u _ \tau } { \partial \nu}
\right) \eta _R \, d\sigma 
 & = &
\int _{ \Omega }  \nabla \eta_R\cdot \left(  | \nabla u_\tau |^2
\alpha   - (2 \alpha \cdot \nabla u _ \tau) \nabla u_\tau \right)\,
dy\\
&\equiv  & I_R.
\end{eqnarray*}
The hypothesis $\epsilon > -1$ is needed to justify integration by parts
with the singular vector field $ \alpha$. 
%\marginpar{Best way to do this? OK}
Now, using the measures defined in (\ref{Ecarleson}),
Proposition \ref{Carleson} and the definition of $\alpha$ we have  
$$
|I_R | \leq C \int _{ \partial \Omega }  f_R(x) \, d\sigma _
 \epsilon 
$$
where $ f_R$ is given by
$$
f_R(x) =  \sup _{ y \in \Gamma (x) , |y |> R}  | \grad u_{\tau} (y) |^2 . 
$$
Notice that for each $x$,  $ \lim _{R \rightarrow \infty }  f_R(x) =0
$. Hence, our 
assumption that $(\nabla u )^*$ is in $L^2 (d\sigma _ \epsilon)$ and the
Lebesgue dominated convergence theorem imply  that $ \lim_{ R
\rightarrow \infty } I_R =0 $. 

We may now let $ \tau \rightarrow 0^+$ and use the dominated
convergence theorem again to obtain the  lemma. 
\end{proof}
%% Our goal is to obtain $L^p$-results for the mixed problem in two dimensions
%%  for a  
%% class of domains that are not included in Sykes \cite{JS:1999} and 
%% Sykes and Brown \cite{SB:2001}. In particular, we will consider some
%% domains for which we cannot solve the mixed problem in $L^2 (d\sigma)$.  
%% As a step in this direction, we will study weighted results for the
%% Neumann and regularity problems in graph  domains.
As a step towards studying the mixed problem in weighted spaces, we 
consider the Neumann and regularity problems in two-dimensions. 
By the \emph{Neumann problem for }$ L^p(w\, d\sigma )$, we mean the
problem of finding a   function $u$ 
which satisfies 
\begin{equation}\label{Neumann}
\left\{ \begin{array}{ll} 
           \Delta u = 0 \qquad & \mbox{in } \Omega \\
           \frac { \partial u } { \partial \nu } = f_N,  \qquad & \mbox{on
           } \partial \Omega  \\
           (\nabla u )^* \in L^ p ( w \, d \sigma )
        \end{array}\right. 
\end{equation}
where, in general, we assume that $f_N$ is taken from  
$L^p( \,w \, d \sigma )$. 
%or from a Hardy space.  
We also study the \emph{regularity problem for }  $L^p(w \, d\sigma)$ where we
look for a function $u$ which satisfies
\begin{equation} \label{Regular}
\left\{ \begin{array}{ll} \Delta u = 0 \qquad & \mbox{in } \Omega \\
u = f_D,  \qquad & \mbox{on
         } \partial \Omega  \\
 (\nabla u )^* \in L^ p ( w \, d \sigma )
         \end{array}\right. 
\end{equation} 
where, in general, we assume that 
 $ d f_D/d\sigma$ is in  
$L^p( w \, d \sigma )$.
%, or from a Hardy space. 

We will obtain  weighted estimates for these problems by applying the
Rellich identity  (Lemma \ref{RellichIdentity}) with a vector field $ \alpha$
which satisfies  $- \alpha \cdot \nu \approx |x|^\epsilon$ for
appropriate  $\epsilon $. Here and in the sequel we will use the
 notation $ A \approx B$ to signify: $c_1A\le B\le c_2A$ for fixed
 constants $0<c_1<c_2<+\infty$. 

\begin{lemma} \label{vecfield} Let $ \Omega$ be a Lipschitz graph 
domain
  with Lipschitz constant $M$ and set $\beta = \arctan (M)>0$. For 
$|\epsilon|<(\pi - 2\beta)/(\pi + 2\beta)$
 there exist $\beta_0 = \beta_0(\epsilon, M)$,
$\beta<\beta_0<\pi/2$,  
%$( \pi/2-|\beta| )/ ( \pi/2+|\beta|)$, 
%there exist a constant $0<c<1$
 and a complex number $a=e^{i\lambda}$ 
 such that the vector field 
$ \alpha(z) =(\Re ( a z^\epsilon), \Im ( a
z^\epsilon ))$ has the following property:  
$$
-|x|^{\epsilon}\le
\alpha (x) \cdot \nu (x)
< - |x|^\epsilon\sin(\beta_0 -\beta)\, , \qquad x \in \partial \Omega. 
$$
\end{lemma}
\comment{
\begin{remark} The vector field $\alpha$ is holomorphic
 (indeed conformal) in a neighborhood of $\overline{\Omega}$.
\end{remark}
}

\begin{proof}
Since $ \beta = \arctan (M)$,
% the  region $ \Omega$ is contained in
%the sector $ \{\, z = r e^{ i\theta } \, :\, r>0,
%\ -\beta < \theta <  \beta + \pi\, \}$.
the outer unit normal $ \nu$ lies in  $\{ (
\cos \varphi , \sin \varphi ) :
 -\pi/2 -\beta \leq \varphi \leq -\pi/2 + \beta \}$.
 Thus, in order to have $ \alpha \cdot \nu < 0$ we need
$\alpha / |\alpha|= (\cos \psi, \sin \psi )$ for some $ \psi $
 in $(\beta, \pi - \beta)$.
To obtain a strictly negative upper bound for $\alpha \cdot \nu$,
it suffices to pick $ \beta _0$ (to be determined later) so that
 $ \beta < \beta _0 < \pi/2$ 
and then
require $\alpha /|\alpha|= ( \cos \psi, \sin \psi )$ for some
$\psi $ in $(\beta_0 , \pi  - \beta_0)$.
% (an interval that is
% strictly contained in $(-|\beta|, \pi + |\beta|)$). 

For $z=re^{i\theta}$ in $\Omega$, we have
 $-\beta<\theta<\pi+\beta$. In order to construct the vector
 field $\alpha$,  given $\epsilon$ as in the hypothesis, we let 
$ \psi ( \theta )$ be a linear function with slope $\epsilon$
 so that $\psi$ 
 maps the interval $[-\beta ,  \pi +  \beta  ] $ onto
$[ \beta_0 , \pi - \beta_0 ]$. 
That is, we let $\psi (\theta) = \epsilon \theta + \lambda$, 
where $|\epsilon| = (\pi - 2 \beta_0 )/(\pi + 2 \beta )$
(we may choose $\psi$ to have either positive or negative slope).
The latter defines $\beta_0$.

We now let 
$\alpha(z) = e^{i\lambda}z^{\epsilon}$, or in polar coordinates, 
$\alpha ( r e^{ i\theta} ) = r^{\epsilon}e^{i\psi (\theta)}$.
%% $$
%% \alpha_{\pm} ( r e^{ i\theta} ) = ( \Re r^{\pm m} e^{ i \ell_{ \pm
%% }( \theta
%% ) } ,  \Im r^{\pm m} e^{ i \ell_{ \pm }( \theta ) }  ) . 
%% $$ 
Then the angle between $ \alpha$ and $ \nu $ will lie in the interval
$( \pi/2 +  \beta_0  -  \beta  , 3\pi /2 - ( \beta_0  -  \beta ))$ and we have 
$$
%\alpha_{\pm} \cdot \nu \leq - \sin ( \beta_0 - \beta ) r^ {\pm m }.
-|x|^{\epsilon} \le \alpha (x)\cdot\nu (x)
 \le - |x|^{\epsilon}\sin (|\beta_0| - |\beta|)\, ,\qquad x \in \partial\Omega.
$$
%We may choose $\psi$ to have either positive or negative slope $\epsilon$.
%In either case, it easy to see that
% $|\epsilon| = (\pi - 2 \beta_0 )/(\pi + 2 \beta )$.
%
%$Concerning the selection of $|\beta_)|$, given 
%$|\epsilon| = (\pi - 2|\beta|)/(\pi + 2|\beta|)$ it suffices to let
%$|\epsilon| = (\pi - 2|\beta_0|)/(\pi + 2|\beta|)$.
\end{proof}

\begin{proposition}[Weighted Rellich Estimates for (\ref{Neumann})
    and (\ref{Regular})]  \label{RellichIneq}  
Let  $\Omega $ be a Lipschitz graph domain with 
Lipschitz constant $ M$ and let $  \beta  =\arctan (M)$. Suppose that
$ |\epsilon |<(\pi - 2 \beta )/(\pi + 2 \beta )$. Then, if
%$ < (  \pi/2 -|\beta| ) / ( \pi/2 + |\beta|)$. Then, if
 $u$ is harmonic in $\Omega$ and  
 $ \nontan{(\nabla u)}  \in L^2 ( \sigma_\epsilon )$, we have 
$$
\int_{ \partial \Omega } \left| \frac { \partial u }{ \partial \nu }
\right | ^2 \,d \sigma _\epsilon \leq
\int_{ \partial \Omega } | \nabla u |^2 \, d\sigma_ \epsilon \leq C 
\int_{ \partial \Omega } \left| \frac { \partial u }{ \partial \nu }
\right | ^2 \,d \sigma _\epsilon
$$
and
$$
\int_{ \partial \Omega } \left| \frac { du }{ d\sigma }
\right | ^2 \,d \sigma _\epsilon\leq 
\int_{ \partial \Omega } | \nabla u |^2 \, d\sigma_ \epsilon \leq C 
\int_{ \partial \Omega } \left| \frac { du }{ d\sigma }
\right | ^2 \,d \sigma _\epsilon . 
$$
where $ C= C(\epsilon, M)$. 
\end{proposition}

\begin{proof} We let $ \alpha $ be the vector field from 
Lemma \ref{vecfield}  and use the Rellich identity from Lemma
\ref{RellichIdentity} 
to obtain  
$$
\int _{ \partial \Omega } | \nabla u |^2 \alpha \cdot \nu - 2 \frac {
\partial u }{ \partial \nu } \alpha \cdot \nabla u\, d\sigma = 0. 
$$
Because we have $ - \alpha \cdot \nu (x)   \approx |x| ^\epsilon$,
standard arguments as in Jerison and Kenig \cite{JK:1982c} yield  the
proposition. 
\end{proof}
%%Corrections to here. 

An important component in establishing uniqueness is the following
local regularity result for solutions with zero data.

\begin{lemma} \label{LocalRegularity}  Suppose that $\Delta u =0 $ in
$\Omega $, $ \frac {\partial u }{\partial \nu } = 0 $ or $ u =0 $ on
$\partial \Omega $ and $ \nontan {(\nabla u )} \in L^1_{loc} (\partial
  \Omega ) $. 
Then, for every bounded set $B\subset \Omega$, we have $ \nabla u \in
L^ 2 ( B)$.  
\end{lemma}

\begin{proof}
We consider the case of Neumann boundary conditions. Dirichlet
boundary conditions may be handled by a similar argument. 

We first observe that since $ \nontan{(\nabla u )} $ 
is in $L^1_{loc}( d\sigma)$, it follows that 
$\nabla u $ is in $L^1(B)$ for any bounded subset  of $ \Omega$, $B$. 

Fix a point $x_0 \in \partial\Omega$ and
let $\eta $ be a smooth cutoff function which is one on $ \ball { x_0}
{2r} $ and zero outside of  $ \ball {x_0} { 4r} $. We let $N(z,y)$ be
the Neumann function for $ \Omega$, thus $N$  is a  fundamental
solution which satisfies homogeneous Neumann boundary
conditions. Dahlberg and Kenig 
\cite{DK:1987} construct the Neumann function on a graph domain by
reflection.  At least formally, we have the representation formula,
\begin{eqnarray}\label{RepFmla}
\eta u (z) & = & \int _{ \Omega \cap \ball {x_0} { 4r}}  N(z,y) ( u (y)
\Delta \eta ( y) +2 \nabla u (y) \cdot \nabla \eta (y) ) \, dy 
\\
& & \quad - \int _{ \sball { x_0 } { 4r} } N  (z,t) u (t)
 \frac {\partial  \eta
  ( t)} {\partial \nu }  \, d\sigma ( t)\, , \qquad z\in \Omega.\nonumber 
\end{eqnarray}

We now show how to estimate each term on the right-hand side
of this formula. Since $\nontan{( \grad u) }( x) $ is in $L^1_{loc} (
\partial \Omega) $, it follows that $u$ is bounded on $\sball {
x_0} { 4r} $. 
As we observed above, we have  $ \nabla u \in L^ 1 (
\Omega \cap \ball { x_ 0 } { 4r})  $, and the Poincar\'e inequality gives 
that $u$ is in $L^1 ( \Omega \cap \ball {x_0} { 3r})$. We also have
that the Neumann function is locally bounded when  $ z \neq y$. Hence,
the integrands in  (\ref{RepFmla}) are in $L^1$ and it is a routine matter to
justify this formula. Since the map 
$ y \rightarrow N(z,y)$ is in the Sobolev space
 $\sobolev { 1 } { 2} (\ball {x_0} r \cap \Omega )$,
 uniformly for $z \in \ball {x_0} {3r}
\setminus \ball { x_ 0 } {2r}$, the estimates for $u $ and $\nabla u$
outlined above together with (\ref{RepFmla}) imply that
 $\nabla u \in L^ 2(\ball {x_0 } r \cap \Omega )  $.
\end{proof}

Next, we recall a classical fact from harmonic function theory, 
a Phragmen-Lindel\"of theorem.  This result is
well-known and the proof is omitted, see e.g. Protter and Weinberger
\cite[Section 9, Theorem 18]{PW:1967}.  
We will need this theorem to control the behavior at infinity of
solutions in our unbounded domains.

In the next result and below, we  will let $ \Omega _\varphi$ denote
the sector  
\begin{equation}\label{Esector}
\Omega _\varphi\, =\, 
\{ re^{i\theta} : r > 0, \  |\theta- \pi/2| < \varphi\}, 
\quad  0 < \varphi< \pi\, .
\end{equation} 
With this normalization, $ \Omega_\varphi$ is a Lipschitz domain with
constant $ M =  |  \tan ( \pi / 2 - \varphi)  | $.
\begin{theorem} [Phragmen-Lindel\"of]  \label{phrag}
 Let $\varphi \in (0, \pi)$. Suppose $v$ is
sub-harmonic in $\Omega _ \varphi$,
 $ v = 0$ on $ \partial \Omega _ \varphi$ and 
$$
v(z) = o( |z|^{ \pi/(2\varphi)})\quad
 \textit{as}\ |z|\to +\infty\, ,\quad z\in \Omega_{\varphi}\, .
$$
Then $ v \leq 0$. 
\end{theorem}

%Here, we continue to use $ \beta = \arctan (M)$.
We now prove uniqueness for the regularity problem
for $L^p(wd\sigma)$, (\ref{Regular}).  
\begin{lemma}\label{RPUniq}
% Suppose that
% $ \Omega \subset \Omega_\varphi$ and let
% $ \epsilon < \pi/(2\varphi) -1/2$. 
Let $\Omega$ be a Lipschitz graph domain with Lipschitz constant
 $M$. Suppose that, for $p\ge 1$,   
 $L^p_{loc}(w\, d\sigma ) \subset L^1_{loc}(d\sigma)$ and that 
that for all surface balls
$\sball x r $ with $r>1$,  we have 
\begin{eqnarray} \label{Uniq2}
\left( \int_{ \sball x r } w(t)^{ -1/ (p-1)} \, d \sigma(t)  \right )^{
\frac{p-1}{p}} \leq C r^{\pi/(2 \beta  + \pi)}, & \textit{ if }p>1;  \\
%\leq C r^ { \frac 1 2 + \epsilon } . 
\left(\inf\limits_{t\in \Delta_r(x)} w(t)\right)^{-1} \le  C r^{\pi/(2
  \beta  + \pi)}, & \textit{ if }p=1
\end{eqnarray}
 where $\beta = \arctan M >0$.
Under these conditions, if $u$ is harmonic in $\Omega$, 
$ \nontan  {(\nabla u )} \in L^ p ( w\,d\sigma )$  and 
$u=0$ on $\partial \Omega $, then $ u = 0$ in $ \Omega$. 
\end{lemma}

\begin{proof}
 Since $  \beta  = \arctan M$  then we have
$ \Omega \subset \Omega_{\tilde\varphi}$, with
 $\tilde\varphi= \beta +\pi/2$.
%%\marginpar{Uniq for $L^1+L^2$?}
 Define $v$ by 
$$
v (x) =\left\{ \begin{array}{ll} |u | , \qquad & \mbox{in } \Omega \\
                                  0, \qquad & \mbox{in } \bar \Omega
				  ^c.
\end{array}\right. 
 $$
The function $v$ is sub-harmonic in all of $ \reals ^2$ and, moreover,
 $v=0$ on $\partial \Omega_{\tilde\varphi}$.  We
verify the growth condition in the Phragmen-Lindel\"of Theorem
 \ref{phrag}.
To do this, suppose that $ z \in \Omega_{\tilde\varphi}$ and set
% $ r= 2 \dtb x$,
$r= |z|$. 
By Lemma \ref{LocalRegularity} we have that $ \nabla u$
is in $ L^2$ of each bounded subset of $ \Omega$, and the same
is true for $ \nabla v$.  
Then, by combining the mean-value property for sub-harmonic functions
with the Poincar\'e inequality,  Proposition \ref{Carleson} (for the
Carleson measure 
$d\mu (y) = \frac{1}{r}\chi_{B_r(z)}\, dy$ with respect to $d\sigma$,
see (\ref{Ecarleson}))    and H\"older's 
inequality, in the case $p>1$,  we obtain 
\begin{eqnarray*}
0\le v(z)  & \leq & \frac 1 { \pi r^2} \int_{\ball z r} v(y) \, dy  \\
        & \leq & \frac C r \int_{ \ball z r } |\grad v (y)|\, dy \\
        &  \leq & C \int _{ \sball x {2 r} } \nontan{  (\nabla u ) } (t) \,
        d\sigma(t)  \\
        & \leq & C \left( \int_{\sball x {2 r} }  \nontan{( \nabla   u )} (t)
        ^p \, w(t) d \sigma (t) \right ) ^{ 1/p} \left( \int_{\sball
        x {2r}}  w^{ -1/ (p-1)} (t) \, d\sigma (t)\right)^{\frac{p-1}{p}}.
\end{eqnarray*}
Here, $x$ is the projection of $z$ onto the boundary, i.e. if $ z=
(x_1, \phi(x_1)+t)$, then $ x = (x_1, \phi(x_1))$. 
Thus, under our assumption (\ref{Uniq2})
it follows that
$$0\le v(z)  \leq C  r ^ {\frac{\pi}{2\tilde\varphi}} 
\le C|z|^{\frac{\pi}{2\tilde\varphi}},\ \ z\in \Omega_{\tilde\varphi}\, ,
$$
as $0\in \partial\Omega_{\tilde\varphi}$. 
%On account of the hypothesis:
%$ \epsilon < \frac {\pi}{ 2\varphi } -\frac 1 2$
We may now apply Phragmen-Lindel\"of Theorem \ref{phrag} and conclude 
that $ v= 0$. The case $p=1$ is treated in a similar fashion.  
%\marginpar{Check on range.}
\end{proof}
\begin{remark}
Using H\"older's inequality, we see that
$$
\int _{ \sball {x} s } |f(t)| \, d\sigma (t)
\leq \left( \int _{ \sball {x} s }|f(t)| ^p \, w(t) \,
d\sigma (t) \right)^{ 1/p}  \left( \int_{\sball {x} s} w(t)^{-1/(p-1)}\,
d\sigma(t)\right)^{\frac{p-1}{p}}.
$$
Thus, we have $ L^p_{loc}( w \, d\sigma)\subset L^1 _{ loc} (d\sigma)$
provided  $w^{-1/(p-1)}$  is in $L^ 1 _{loc}(d\sigma)$. It is easy to
see that this will hold for the weight $ w(t) = |t|^\epsilon$ if $
\epsilon < p-1$. 
\end{remark}

Finally, we give uniqueness for the weighted Neumann problem
 for $L^p(w\,d\sigma)$, (\ref{Neumann}). 

\begin{lemma} \label{NPUniq} 
Let $\Omega$ be a Lipschitz graph domain with Lipschitz constant $M$.
Suppose $w$ satisfies the hypotheses of
 Lemma \ref{RPUniq}. If $u$ is a solution of the Neumann problem with 
$\nontan{(\nabla u)} \in L^p(w\, d\sigma)$, $p\ge  1$, and
 $\frac{\partial u}{\partial\nu}=0$ a.e. $\partial\Omega$,
 then $u $ is constant. 
\end{lemma}

\begin{proof} We consider $v$, the conjugate harmonic function to
$u$. Since $u$ has normal derivative zero at the boundary, by
the Cauchy-Riemann equations   we have that $v$ is 
constant on the boundary, and we may assume this constant is 0. Now,
$\nontan{(\nabla v)} \in L^p(w\, d\sigma)$ (since this is the case for $u$)
so that $v$ satisfies the hypotheses of Lemma \ref{RPUniq} 
and hence it is zero. But this implies that $u$ is constant.  
\end{proof}

\section{The $L^2$-Neumann and regularity problems with 
 power weights}
\setcounter{equation}{0} 

In this section we prove existence of solutions for  the 
Neumann problem and the regularity problem for $L^2(w\, d\sigma)$
when the weight is a power of $|x|$. 
A different proof  of this result was given by Shen \cite{ZS:2005} for
bounded domains in dimension three and higher; it is likely that
Shen's argument can be adapted to graph domains in two
dimensions. Here we give a  proof based on the Rellich estimates
(Proposition \ref{RellichIneq}).  

%%%to here. 

Before continuing, we record a few basic facts about the power weights
$ |x|^\epsilon$ and their relation to the Muckenhoupt class
$A_p(d\sigma)$.  A weight (that is a non-negative measurable and locally
 integrable function)
$w$ is a member of the class $A_p(d\sigma)$ if and only if for all
 surface balls
%%arcs 
$ \Delta\subset \partial \Omega$, we have
$$
\int _\Delta w \, d \sigma \left( \int_\Delta w^{\frac{-1}{p-1}}\, d\sigma
  \right)^{p-1}
\leq C.
$$
The best constant $C$ in this inequality is called the $A_p$-constant
for $w$. The following simple lemma tells us that  $|x|^\epsilon $ is in
 $A_p(d\sigma)$ if and only if  $-1<\epsilon < p-1$, see
\cite[p.~218]{ES:1993}. 

\begin{lemma} \label{measure-property} \label{51}
 For  $ \epsilon > -1$,  the boundary measure 
  $\sigma_\epsilon$ (see (\ref{Ecarleson})) satisfies
$$
\sigma_\epsilon ( \Delta _r (x)) \approx  r \max
(|x|,r )^\epsilon, \ x\in \partial\Omega. 
$$
\end{lemma}
The proof is omitted. Next, we  establish existence of solutions in 
 sectors, see (\ref{Esector}). 
%namely 
%\begin{equation}\label{Esector}
%\Omega_\varphi = \{ (x_1 , x_2) : x_ 1+ i x_2 = r e^{ i \theta } , \ |
%\theta - \pi/2 | < \varphi \}.
%\end{equation}
%With this normalization, $ \Omega_\varphi$ is a Lipschitz domain with
%constant $ M =  |  \tan ( \pi / 2 - \varphi)  | $. 
%We describe when we can solve the mixed problem in a
%sector. 

\begin{proposition} Let $ \epsilon$ satisfy $ 1-\frac \pi { 2\varphi} <
\epsilon < 1$, and let $\sigma_{\epsilon}$ be as in (\ref{Ecarleson}).
 Then, we may solve the Neumann problem (\ref{Neumann}) and 
 the regularity problem (\ref{Regular})
for $L^2 (d \sigma _\epsilon )$ in the sector $\Omega_\varphi$, and we
have
\begin{eqnarray*} 
\| ( \nabla u )^*\| _{ L^2 (\partial\Omega_{\varphi}, d\sigma_ \epsilon)}
&\leq C\|f_N\|_{ L^2(\partial\Omega_{\varphi}, d\sigma_\epsilon)}
&\quad \textit{(for \ (\ref{Neumann}))}\, ,
\\
\\
\| ( \nabla u )^*\| _{ L^2 (\partial\Omega_{\varphi}, d\sigma_ \epsilon)}
&\leq C\|f_D\|_{ L^2(\partial\Omega_{\varphi}, d\sigma_\epsilon)}
&\quad \textit{(for \ (\ref{Regular}))}\, .
\end{eqnarray*}
\label{confsectors}
\end{proposition}

\begin{proof} The proof goes along the same lines as the proof 
 of Proposition \ref{sectors}, which deals with 
 the mixed problem in sectors, with the role of Theorem
 \ref{unweightedmp} being played by the classical results on 
the regularity and Neumann problems for the upper half plane. 
 We omit the argument here
in order to save space.
\end{proof}
%\typeout{Use conformal map.}
%\marginpar{Use conformal map?}

The main result of this section is the following theorem. 
\begin{theorem} \label{39.3} Let $ \Omega $ be a Lipschitz graph domain with 
$ \| \phi' \| _\infty = \tan  \beta  $, $0 <  \beta  < \pi/2$.
 Then, for
  $\epsilon$ satisfying $|\epsilon|< (\pi - 2 \beta )/( \pi + 2 \beta )$
 and for $\sigma_{\epsilon}$ as in (\ref{Ecarleson}),
 we have that the
Neumann problem (\ref{Neumann}) and the regularity problem (\ref{Regular}) 
 for $L^2(d\sigma_{\epsilon})$ are uniquely solvable.
\end{theorem}

\begin{proof} We apply the method of layer potentials 
as in Verchota \cite{GV:1984}.
 Note that for
$ \epsilon \in ( -1,1)$, $|x|^\epsilon $ is an $A_2(d\sigma)$ weight (see
Lemma \ref{measure-property}) so that singular integrals are bounded
on $L^2(\partial\Omega, d\sigma_\epsilon)$. Applying
% the Rellich identity Lemma
%\ref{RellichIdentity} with vector field $\alpha$ as in
 Lemma \ref{vecfield} and Proposition \ref{RellichIneq}
 leads to the estimates
\begin{equation}
\int _{ \partial \Omega} \left| \frac { \partial u } { \partial \nu } \right|^2
\, d\sigma_\epsilon \approx
\int _{ \partial \Omega} \left | \frac { d u } { d\sigma  }\right|^2
\, d\sigma_\epsilon\, , \label{Rellich-for-NR}
\end{equation}
which hold for $u$ harmonic in $ \Omega$ with  $( \nabla u ) ^* \in L^2 (
\sigma_\epsilon)$ and for $ |\epsilon|<(\pi - 2 \beta )/(\pi + 2 \beta )$.
 By Proposition \ref{confsectors}, we can
solve the Neumann and regularity problems for 
 $ L^2 (\partial\Omega_{\varphi}, d\sigma_\epsilon)$, where $\epsilon$
 is as in the hypothesis and
 $0<\varphi\leq ( \beta  + \pi/2)/2$ (for which we have
 $1-\pi/2\varphi \leq (2 \beta -\pi)/(\pi + 2 \beta )<\epsilon$).  Estimate 
(\ref{Rellich-for-NR}) and the method of continuity
(see Brown \cite[Lemma 1.16]{RB:1994b} and Gilbarg and Trudinger
\cite[Theorem 5.2]{GT:1983})  now lead to the existence of a solution
in general graph domains.

Uniqueness follows from Lemma \ref{NPUniq} (or Lemma 
\ref{RPUniq}). It is easy to check that
%solutions with $\nontan{(\nabla u)} \in L^2 ( \sigma_\epsilon)$
$w(x) = |x|^{\epsilon}$ for $\epsilon$ as in the hypothesis 
(indeed, $\epsilon >(2 \beta -\pi)/(\pi + 2 \beta )$),
satisfies the growth condition needed in these Lemmata.
\end{proof}

\section{The mixed problem with power weights.}
\setcounter{equation}{0}
\label{section-l2-weight}

The main result of this section is Theorem \ref{MixedProb2} where we
prove existence and uniqueness for the solutions of 
the weighted mixed problem in $L^2$ on a Lipschitz graph domain.  
There will be a restriction on the size of the Lipschitz constant. 
We first discuss the unweighted case. 

\begin{lemma}[Rellich Estimates for the unweighted mixed problem]\label{REUMP}
%Let $\delta >0$ and
 Let 
$\Omega =\{x_2>\phi (x_1)\}$, $N$, $D$ be a standard domain for the
 mixed problem,
see (\ref{EDN}). Suppose that $\phi_{x_1} >\delta>0$ on $N$ and
$\phi_{x_1} <-\delta<0$ on $D$.

Then, if $u$ is harmonic in $\Omega$ and
 $\nontan{(\nabla u)} \in L^2(\partial\Omega,d\sigma)$, we have
\begin{equation}\label{EREUMP}
\int_D \left(\frac {du}{d\sigma} \right) ^2 \, d \sigma + \int_N
\left ( \frac {\partial u }{\partial \nu} \right)^2 \, d \sigma \approx
\int_N \left (\frac {d u }{d\sigma}\right)^2 \, d \sigma + \int_D \left(
\frac {\partial u }{ \partial \nu } \right)^2 \, d \sigma.
\end{equation}  
\end{lemma}
\begin{proof}
The proof follows from 
the Rellich identity with vector field $ \alpha = e_1$,
see Brown \cite[Lemma 1.7]{RB:1995a}. 
\end{proof}
\begin{theorem} \label{unweightedmp}
Let $ \delta >0$ and $ \Omega$, $N$, $D$ be a standard Lipschitz
graph domain for the mixed problem, see (\ref{EDN}). Suppose that $
\phi_{x_1} > \delta>0 $ on $N$ and $ \phi_{x_1} < -\delta<0 $ on $D$. 
Then, if $g_N \in L^2(N,d\sigma)$ and $d g_D/d\sigma \in L^2(D,d\sigma)$,
there exists a unique solution of the mixed problem (\ref{Mixed}) for
$L^2(d\sigma)$ in $\Omega$. 
% (\ref{Mixed}) for $L^2(
%\sigma)$ 
\comment{
\begin{equation}\label{EunwMP}
 \left\{
\begin{array}{ll}
\Delta v = 0 \qquad & \mbox{in } \Omega \\
v= g_D\qquad & \mbox{on } D \\
\frac { \partial v }{ \partial \nu } = g_N \qquad & \mbox{on } N \\
(\nabla v )^* \in   L^2( d\sigma ).
\end{array}
\right.
\end{equation}
}
Moreover, $v$ satisfies 
\begin{equation}\label{Egradv}
\| ( \nabla v )^*\| _{ L^2 (  d\sigma)} 
\leq C
\left ( \| \frac { dg _D}{d\sigma} \| _{L^2 ( D,d\sigma)} + \| g_N \|
_{L^2 (N,d\sigma )} \right) 
\end{equation}
\end{theorem}

\begin{remark}
In particular, we have that the unweighted mixed problem is uniquely
 solvable (with estimates for $(\nabla v )^*$) in all convex sectors
$\Omega_{\psi}$, $0<\psi<\pi/2$, see (\ref{Esector}). 
\end{remark}
%39.5
\begin{proof}
This is an extension of the results of Brown \cite{RB:1994b} to two
dimensions and to Lipschitz graph domains. The first step is to
observe that the equivalence (\ref{EREUMP}) quickly leads to uniqueness. 
%%\marginpar{Check powers}
For if $u$ is a solution of  (\ref{Mixed}) with $u=0$ on $D$
and $ \frac
{\partial u }{ \partial \nu} =0$ on $N$, then 
(\ref{EREUMP}) implies $u=0$ on $\partial \Omega$. 
 Lemma \ref{RPUniq} now yields $u=0$ in
$\Omega$. 

To prove existence, we find solutions in a sector. This is
easy by symmetry and can be done as in Brown \cite[Lemma
1.16]{RB:1994b}. Then we may use   
 the method of continuity and again (\ref{EREUMP}) to
establish solutions in more general Lipschitz graph domains. 
\end{proof}

We now turn our attention to the weighted mixed problem in
 sectors that are not necessarily convex. 
%This is a step towards the general weighted problem, just as
%in the previous proposition
%we used the unweighted mixed problem for sectors as a
%step towards the solution of the unweighted 
% mixed problem in more general domains. 
For the boundary of a sector $\Omega_\varphi$ as in (\ref{Esector}), we
write: $\partial\Omega_\varphi = D_\varphi \cup N_\varphi$, with 
$D$ and $N$ as in (\ref{EDN}).     
%39.9

\begin{proposition} \label{sectors}
Let $ \varphi \in ( \pi/2, \pi)$ and suppose that
$1-\frac {\pi}{2\varphi} < \epsilon < 1$. Let $\sigma_\epsilon$ be
as in (\ref{Ecarleson}).  
% and $ \epsilon > -1$.
Then, if $f_N \in L^2(N_{\varphi}, d\sigma_{\epsilon})$
and $d f_D/d\sigma \in L^2(D_{\varphi}, d\sigma_{\epsilon})$,
there exists  a unique solution $u$ of the mixed problem 
%the weighted mixed problem in $\Omega_\varphi$ 
%for
%$L^2(\sigma_\epsilon)$  in $ \Omega_\alpha$ 
\begin{equation}\label{EwMP}
 \left\{
\begin{array}{ll}
\Delta u = 0 \qquad & \mbox{in } \Omega_\varphi \\
u= f_D\qquad & \mbox{on } D_\varphi \\
\frac { \partial u }{ \partial \nu } 
= f_N \qquad& \mbox{on } N_\varphi \\
(\nabla u )^* \in   L^2(\partial\Omega_\varphi, d\sigma_\epsilon ).
\end{array}
\right.
\end{equation}
Moreover, $u$ satisfies
\begin{equation}\label{Eestimate} 
\int_{\partial\Omega_\varphi} (\nabla u )^* (x)^2 \, d\sigma_\epsilon
 \leq C \left( \int_{N_\varphi} \left |f_N(x)\right|^2 \, d \sigma_\epsilon
  +\int_{D_\varphi} \left |\frac {df_D}{d\sigma}(x)\right |^2 \,
 d \sigma_\epsilon\right)
\end{equation}
\end{proposition}

\begin{proof} 
We use a conformal map to reduce the weighted Mixed problem (\ref{EwMP}) to
an unweighted Mixed problem on a convex sector, then we apply Theorem
\ref{unweightedmp}.
 
Let $s=1-\epsilon$ with $\epsilon$ as in the hypothesis, so that we have 
\begin{equation}\label{Es}
0<s<\pi/2\varphi<1.
\end{equation} Define $ \phi_s(z) $ to be the  conformal map
$$
\eta=\phi_s(z) =i (-iz)^s.  
$$
Thus, $ \phi_s $  maps $\Omega _\varphi$ to $\Omega_{s\varphi}$,
$D_\varphi$ to $D_{s\varphi}$ and $N_\varphi$ to $N_{s\varphi}$.
Note that if we let $\partial\phi_s$ denote
 the complex derivative of $\phi_s$, that is
\begin{equation}\label{ECderiv}
\partial\phi_s (z) = \frac{1}{2}\left(\frac{\partial\phi_s}{\partial x_1}
 + \frac{1}{i}\,\frac{\partial\phi_s}{\partial x_2}\right)
\end{equation}
we have
\begin{equation}\label{Edphi}
\frac{1}{|\partial\phi_s|}d\sigma= \frac 1 s d\sigma_\epsilon\, . 
\end{equation} 
On account of (\ref{Es}) it follows that Theorem \ref{unweightedmp}
applies to $\Omega_{s\varphi}$.
 
Given $f_D$ and $f_N$ as in (\ref{EwMP}), we define $g_D$ and $g_N$
on $\partial\Omega_{s\varphi}$ as follows:
\begin{equation}\label{Eg}
g_N(\eta) = \frac{f_N}{|\partial\phi_s|}\circ\phi_s^{-1}(\eta)\, ;\quad
g_D(\eta) = f_D\circ\phi_s^{-1}(\eta).
\end{equation}  
We verify that $g_N$ and $g_D$ satisfy the hypotheses of Theorem
\ref{unweightedmp}. Indeed, it is immediate to see that
\begin{equation}\label{EgN} 
\int_{N_{s\varphi}}|g_N|^2d\sigma =\frac 1 s
\int_{N_{\varphi}}|f_N|^2d\sigma_{\epsilon}\, .
\end{equation}
Similarly, we have
\begin{equation}\label{EgD1}
\int_{N_{s\varphi}}\left|\frac{dg_D}{d\sigma}\right|^2d\sigma =
\frac 1 s
\int_{N_{\varphi}}\left|\frac{df_D}{d\sigma}\right|^2d\sigma_{\epsilon}\,. 
\end{equation}
%It follows
%$\frac{dg_D}{d\sigma} \in L^1_{loc}(\partial\Omega_{s\varphi}, d\sigma)$,
% and by the fundamental theorem of calculus we have that $g_D$ is 
%continuous in $D_{s\varphi}$. In particular, we have: 
%$$g_D \in L^2_{loc}(\partial\Omega_{s\varphi}, d\sigma).$$
%We are left to show
%$$I_{R_o} = \int\limits_{D_{s\varphi}\, \cap \{\, |\eta|\,>R_o\}}
%\!\!\!\!\!\!\!\! |g_D|^2\,d\sigma
%\ <+\infty.$$
%
% Changing variable: $\eta = \phi_s(z)$, we obtain
%
%$$I_{R_o} = \int\limits_{D_{s\varphi}\, \cap \{\, |z|\,>R_o^{1/s}\}}
%\!\!\!\!\!\!\!\! |f_D|^2|\partial\phi_s|^2\,d\sigma_\epsilon\, .$$
%Now, $|\partial\phi_s|^2(z) \approx |z|^{-2\epsilon}$, thus
%$$I_{R_o} \le C||f_D||_{L^2(D_\varphi, d\sigma_\epsilon)}$$
%and it follows
%\begin{equation}\label{EgD2}
%g_D \in L^2(D_{s\varphi}, d\sigma).
%\end{equation}
 
By Theorem \ref{unweightedmp} the unweighted Mixed problem
 on $\Omega_{s\varphi}$ has a unique solution $v$ which satisfies
(\ref{Egradv}). 
We now pull this solution back to $\Omega_\varphi$ by  defining
$u=v\circ \phi_s$. 
Then, $u$ is harmonic in  $ \Omega_\varphi$ and satisfies:
$$
u= f_D\quad \mathrm{on}\ D_\varphi\, ; \qquad 
\frac{\partial u}{\partial \nu} = f_N\quad \mathrm{on}\ N_\varphi\, .
$$ 
Moreover, on account of (\ref{Egradv}), (\ref{EgN}) and 
(\ref{EgD1}) we have  
%$$
%\int_{N_{s\varphi} }\left| \frac {\partial v }{ \partial \nu}
%  \right|^2\, d\sigma =
%\int_{N_{\varphi}}\left| \frac {\partial u  }{\partial \nu }\right|^2
%  d\sigma_\epsilon 
%$$
%$$
%\int_{D_{s\varphi} }\left| \frac {d v }{ ds}
%  \right|^2\, d\sigma =
%\int_{D_{\varphi}}\left| \frac {d u  }{ds  }\right|^2
%  d\sigma_\epsilon 
%$$
%$$
%\int_{\partial \Omega_{s\varphi} }|\nabla v|^2 
%\, d\sigma =
%\int_{\partial\Omega _{\varphi}}|\nabla u|^2 
%  d\sigma_\epsilon ,
%$$
%%Since we have $ \frac 1 {|\partial \phi_\beta (x)|} =
%% |\beta||x|^{1-\beta} $, the
%%measure $ \frac 1 {|\partial \phi_\beta|} \, d\sigma $ is
%%comparable to $ d\sigma_{1-\beta}$. If we require that $ \alpha \beta
%%< \pi/2$, then we know that a solution of the unweighted Mixed
%%problem exists from Theorem \ref{unweightedmp}. 
%%Thus, given data for the $L^2
%%(\sigma_{1-\beta} )$-Mixed problem, we may map to $
%%\Omega_{\alpha\beta}$, find a solution and then map back to find a
%%function $u$ which is a solution to  Mixed problem for $L^2 (
%%\sigma_{1-\beta})$ in
%%$ \Omega_\alpha$. The condition that $ \alpha \beta < \pi/2$ which is
%%needed in Theorem \ref{unweightedmp} implies that we 
%%%obtain a solution to the weighted Mixed problem for $ 1-\beta >
%%1-\pi/(2\alpha)$. 
%and we obtain
\begin{equation}\label{Eincomplete}
% \int_{\partial\Omega_{s\varphi}} (\nabla v )^* (x)^2 \, d\sigma
\|\nontan{(\nabla v )}\|_{L^2(\partial\Omega_{s\varphi}\, ,d\sigma)}
 \leq C  \int_{N_\varphi} \left |f_N(x)\right|^2 \, d \sigma_\epsilon
  +\int_{D_\varphi} \left |\frac {df_D}{ds}(x)\right |^2 \, d \sigma_\epsilon.
\end{equation}

Now we consider non-tangential maximal function estimates for $\nabla
u$. We
apply the Cauchy integral formula to the complex derivative
of $u$ (which is analytic in $\Omega_{\varphi}$) and obtain: 
\begin{equation}\label{ECauchy}
\partial u (z)\, =\,
\frac{1}{2\pi i}\int _{ \partial \Omega_\varphi} \frac {\partial u (\zeta)}
{ z-\zeta} \, d\zeta \, =\,
\frac{1}{2\pi i}
\int _{ \partial \Omega_\varphi} \frac {\partial \phi_s (\zeta)
 \partial v
  (\phi_s (\zeta))}{ z-\zeta} \, d\zeta \, ,\qquad  
 z\in \Omega_\varphi\, .
\end{equation}
%We let $\xi =\phi_s (\zeta)$ and  rewrite this integral as 
%$$
%\int _{ \partial \Omega _{s\varphi}} \frac { \partial v
%  (\xi)}{z- \phi_s^{-1}(\xi)} \, d\xi .
%$$
%In the domain $ \Omega_{\alpha \beta}$, we have non-tangential maximal
%function estimates. We may use these to evaluate the last integral using
%the residue theorem and obtain that 
%$$
%\frac {1} {2\pi i } \int _{ \partial \Omega _\alpha}
% \frac {\partial \phi_\beta (y) \partial u
%  (\phi_\beta (y) )}{ x-y} \, dy = \partial \phi_\beta (x)\partial   u
%(\phi_\beta (x)) =\partial v (x)  .
%$$
%A similar representation holds for $ \bar \partial u$.  
%From this representation 
%and the theorem on the Cauchy integral
%\cite{CMM:1982}, 
It follows
\begin{equation}\label{ECauchynontan}
\nontan{(\nabla u)}(x) = 2\nontan{(\partial u)}(x) \le
\nontan{\bigg(K_{\varphi}(\partial\phi_s\cdot(\partial v)\circ\phi_s)\bigg)}(x)\, ,
\quad \mathrm{a\,.\ e.\ } \, x\in \partial\Omega_{\varphi}\, ,
\end{equation} 
where $K_{\varphi}$ denotes the Cauchy integral on $\Omega_{\varphi}$.
By the theorem of Coifman, McIntosh and Meyer \cite{CMM:1982} on the  boundedness of the Cauchy integral  we have
\begin{equation}\label{ECauchyweight}
\|\nontan{(K_{\varphi}\psi)}\|_{L^2(\partial\Omega_{\varphi},d\sigma_{\epsilon})}
\le C \|\psi\|_{L^2(\partial\Omega_{\varphi},d\sigma_{\epsilon})}.
\end{equation}
This uses that $d\sigma_\epsilon$ is in $A_2(d\sigma)$. 
 Combining these last two inequalities we obtain
\begin{equation}
\|\nontan{(\nabla u)}\|_{L^2(\partial\Omega_{\varphi},d\sigma_{\epsilon})}
\le C \|\partial\phi_s\cdot(\partial v)\circ\phi_s)\|_{L^2(\partial
\Omega_{\varphi},d\sigma_{\epsilon})}\, =\,
C\|\partial v\|_{L^2(\partial\Omega_{s\varphi},d\sigma)}\, , 
\end{equation}
where the last equality was obtained by performing
 the change of variables $\phi_s(\zeta) = \eta$, see also 
(\ref{Edphi}). 
% where $w(\zeta) = |\partial \phi_s (\zeta)| \in A_2(\partial\Omega_\varphi)$.
%
%The latter may be re-written as
%$$
%(\nabla u )^ *(x) \leq K(|\nabla v|\circ \phi_s\, |\partial \phi_s|)(x)\, ,
%$$
%where $K$ denotes the Cauchy transform on the Lipschitz curve
%$\partial\Omega_\varphi$.
%In view of the weighted $L^2$-regularity of the Cauchy integral, see [??],
%we have
%\begin{eqnarray*}
%\|(\nabla u)^*\|_{L^2(\partial\Omega_\varphi, d\sigma_\epsilon)}^2
% & \le C \int\limits_{\partial\Omega_\varphi} |\nabla v|^2 (\phi_s(t))
%|\partial\phi_s|^2 d\sigma_\epsilon \\ 
% & = C \int\limits_{\partial\Omega_\varphi} |\nabla v|^2 (\phi_s(t))
%|\partial\phi_s| d\sigma \\ & =
% C \|\nabla v\|_{L^2(\partial\Omega_{s\varphi}, d\sigma)}^2\, .
%\end{eqnarray*}
%C \int\limits_{\partial\Omega_{s\varphi}} |\nabla v|^2 d\sigma^2  
This, together with (\ref{Eincomplete}) yields (\ref{Eestimate}). 

The argument is reversible, so we may conclude uniqueness in
$\Omega_\varphi$ from uniqueness in $ \Omega_{ s\varphi}$. 
\end{proof}

\begin{remark} It is well-known that non-tangential maximal function
estimates for harmonic functions behave nicely under conformal
mapping, see  Kenig \cite{CK:1980} and Jerison and Kenig
\cite{JK:1982d}. The previous Lemma, however, considers the non-tangential
maximal function of the gradient. The estimates in this case  appear
to be more involved.  
\end{remark}

We now construct holomorphic vector fields which allow us to use the
Rellich identity of Lemma \ref{RellichIdentity} to obtain
 Rellich estimates for the weighted mixed problem. 

\begin{lemma}\label{LVFMP}
 Suppose $ \Omega= \{ x_2 >
  \phi (x_1)\}$, $N, D$ is a standard Lipschitz graph
  domain for the mixed problem, with Lipschitz constant $M$.
Let $\beta =\arctan M >0$. Assume $\beta<\pi/4$.

Then, for $2\beta/(\pi - 2\beta)<\epsilon<1$ there exist
 $\beta_0 =\beta_0(\epsilon, M)$,
 $\beta<\beta_0<(\pi -2\beta)/2$, and a complex number 
$a=e^{i\lambda}$ such that the vector field 
$\alpha(z) = (Re(az^{\epsilon}), Im(az^{\epsilon}))$ satisfies
\begin{equation}\label{VF*1}
-|x|^{\epsilon}\leq \alpha (x)\cdot \nu (x)  <  -|x|^ \epsilon \sin( \beta_0 -\beta)\, ,  \qquad
 x\in N\, ;
\end{equation}
\begin{equation}\label{VF*2} 
|x|^{\epsilon}\geq \alpha (x) \cdot \nu (x)  > |x|^ {\epsilon} \sin( \beta_0 - \beta)\, , \qquad
 x\in D. 
\end{equation}
\end{lemma}
\begin{proof}
The proof goes along the same lines as Lemma \ref{vecfield}. 
The outer unit normal $\nu$ lies in $\{ (\cos\varphi, \sin\varphi)\, :\ 
-\pi/2-\beta \le\varphi\le-\pi/2+\beta\}$. On account of
(\ref{EDN}) and (\ref{Lipconst}) we have that $N$ is contained
 in the sector $\{x=re^{i\theta}\,:\, -\beta<\theta<\beta\}$,
 whereas $D$ is contained in
 $\{x=re^{i\theta}\,:\, \pi-\beta<\theta<\pi+\beta\}$.
Thus, in order to have $\alpha\cdot\nu<0$ on $N$ we need
$\alpha/|\alpha| = (\cos\psi, \sin\psi)$ for some 
$\psi \in (-2\pi+\beta, -\pi-\beta)$, whereas
$\alpha\cdot\nu>0$ on $D$ requires
$\psi \in (-\pi+\beta, -\beta)$. To obtain a strictly
 negative upper bound for $\alpha\cdot\nu$ on $N$ 
 and a
 strictly positive lower bound for
 $\alpha\cdot\nu$ on $D$, we pick $\beta_0$ (to be
 selected later) so that $\beta<\beta_0<\pi/2$ and
 then, with the same notations as above, require that,
 for $x=re^{i\theta} \in N$, $\psi(\theta)$ lie in
$[\beta_0, \pi-\beta_0]$ whereas, for
 $x\in D$, we require that $\psi(\theta)$ lie in
 $[\pi+\beta_0, 2\pi-\beta_0]$.
 To this end, given $\epsilon$ as in the hypothesis, we let
 let $\psi(\theta)$ be a linear function with slope
$\epsilon$ and choose $\beta_0$ so that $\psi$ maps
the point $\beta$ to $\pi-\beta_0$ and the point
$\pi -\beta$ to $\pi +\beta_0$. Writing
 $\psi(\theta)=\epsilon\theta + \lambda$, we define
$\alpha (z) = e^{i\lambda}z^{\epsilon}$, that is
$$\alpha(re^{i\theta}) = r^{\epsilon}e^{i\psi (\theta)}\, $$
where
$\lambda= \pi-\pi\beta_0/(\pi-2\beta)$
 and $\epsilon = 2\beta_0/(\pi - 2\beta)$ (note that the latter
 defines $\beta_0$).
This construction, however, grants that $\alpha\cdot\nu$
has the desired sign only near the endpoints 
$\theta = \beta$ (for $N$) and $\theta = \pi -\beta$ 
 (for $D$). In order to make sure that $\alpha\cdot\nu$
keeps the desired sign all the way through the two other
 endpoints we need to restrict the range for the selection
of $\beta_0$ to:
$$\beta<\beta_0<\frac{\pi}{2}-\beta\, .$$
Then the angle between $\alpha(x)$ and $\nu(x)$ will lie in the intervals
$(-3\pi/2+(\beta_0-\beta), -\pi/2-(\beta_0-\beta))$, for $x\in N$,
 and in  $(-\pi/2+(\beta_0-\beta), \pi/2-(\beta_0-\beta))$, for $x\in D$,
 so that (\ref{VF*1}) and (\ref{VF*2}) hold. 
\end{proof}

\begin{proposition}[Weighted Rellich estimates for the  mixed problem]\label{rellich} Let  $\Omega$, $N$, $D$ be a standard Lipschitz 
graph domain for the mixed problem (\ref{Mixed}), with
 Lipschitz constant $M<1$. Let $ \beta= \arctan
M$. Then, for $\sigma_\epsilon$ as in (\ref{Ecarleson}), for
 $ \epsilon$
in the range $ 2\beta/( \pi - 2 \beta) < \epsilon < 1$
 and for $u$ harmonic with 
$( \nabla u )^* \in L^2 ( \sigma_ \epsilon)$, 
we have 
$$
\int_{ \partial \Omega} | \nabla u |^2 \, d \sigma_\epsilon
\leq C \left[ \int_N \left( \frac { \partial u }{\partial \nu } \right
  ) ^2 \, d \sigma_ \epsilon
+ \int _D \left ( \frac { du }{ d\sigma}\right) ^2 \, d\sigma_\epsilon
\right].
$$
\end{proposition}

\begin{proof} We use the identity of Lemma \ref{RellichIdentity} with
  the vector 
  field constructed in Lemma \ref{LVFMP} and standard manipulations
  involving the boundary terms as in Brown \cite[Lemma 1.7]{RB:1994b}.
 The key point is
  that since $\alpha\cdot \nu$ changes sign as we pass from $N$ to
  $D$, we can estimate the full gradient of $u$ on the boundary by the
  data for the mixed problem. 
%\typeout{This is very sketchy.}
\end{proof}

\begin{theorem}\label{MixedProb2} Suppose $\Omega$, $N$ and $D$ is a
standard Lipschitz graph domain 
for the mixed problem, with Lipschitz constant $M<1$. Let $ \beta =
\arctan(M)$. Given $ \epsilon$ which satisfies
$$
 2 \beta/( \pi -2 \beta) <\epsilon <1, 
$$
there is a unique solution $u$ to the $L^2 ( \sigma_\epsilon)$ mixed
problem (\ref{Mixed}).
% in $\Omega$, $N$ and $D$ and this solution
Moreover, $u$ satisfies
$$
\int_{\partial \Omega} ( \nabla u )^*(x) ^2\, d\sigma _ \epsilon
\leq
 \left( \int_D \left| \frac {df_D}{d\sigma} \right| ^ 2 \, d \sigma_
\epsilon
+\int_N |f_N|^2 \, d\sigma
_ \epsilon \right).
$$
\end{theorem}

\begin{proof}
For $\epsilon$ in the given range, Lemma \ref{LVFMP} and Proposition 
\ref{rellich}
imply that
$$
\int_{ \partial \Omega} |\nabla u | ^2 \, d\sigma_\epsilon
\leq C \left( \int_D \left| \frac {du}{d\sigma} \right| ^ 2 \, d \sigma_
\epsilon
+\int_N \left|\frac { \partial u} {\partial \nu } \right|^2 \, d\sigma
_ \epsilon \right).
$$
Combining
Proposition \ref{sectors} which gives existence in sectors
with  the method of continuity (Brown \cite{RB:1994b}, Gilbarg and
Trudinger \cite{GT:1983}), we obtain the existence of a solution to
the mixed problem with data in $L^2 ( d\sigma_\epsilon)$.
%%\marginpar{Uniqueness?}
%%Do we need more detail--see notes?

Next, we consider uniqueness. If $u$ is a solution of the mixed problem
 with zero data, then the Rellich estimates (\ref{EREUMP}) imply
 that $u$ is a solution of the 
regularity and Neumann problems for $L^2(d\sigma_{\epsilon})$ with
 zero data. By the results in Section 2 it follows that $u$ is zero.
\end{proof}

\section{The   regularity and Neumann problems in
  $H^1(d\sigma_\epsilon)$.} 
\setcounter{equation}{0}

In this section we assert the existence of solutions for the Neumann
 problem when the data is in $H^1 ( d\sigma_ \epsilon)$, and for
 the regularity problem when the data has one derivative in
 $H^1 ( d\sigma_ \epsilon)$.
 This follows the work of Dahlberg and Kenig \cite{DK:1987}; 
%and is less
%involved than what follows for the mixed problem.
thus, we shall be brief. 
We first recall the definition of the Hardy
spaces $H^1(d\sigma_ \epsilon)$ and $H^{1,1}(d\sigma_ \epsilon)$. 

Let $\epsilon>-1$.
We say that $a$ is an $H^1(d\sigma_\epsilon)$-{\em atom for} 
% for the mixed problem
$ \partial \Omega$ if $a $ is supported on a surface ball $
\Delta_s(x)$, $\int a\, d\sigma =0$  and  $\|a\|_\infty \leq \sigma _
\epsilon ( \Delta _s (x) ) ^{-1}$.
We remark that for $\epsilon \leq 0$,  $L^1(d\sigma_{\e}) \subset
L^1_{loc}(d\sigma)$. Thus we define the space $H^1 ( d\sigma_
\epsilon)$, for $\e\le 0$, to be  the collection of 
functions that are 
%nts in $(Lip_c(R^2))'$ (the dual space of
% of Lipschitz continuous functions with compact support) that can be
represented as 
\begin{equation}\label{EH1}
f(x) = \sum _{ j=1}^ \infty \lambda _j a _ j (x) 
\end{equation}
where  
%convergence is meant in distribution sense,
$\{a_j\}_{j\in\mathbf N}$ is a sequence of $H^1(d\sigma_\epsilon)$-atoms
for $\partial\Omega$ and  the coefficients $ \{\lambda _j\} $ 
satisfy $\sum_{ j=1}^ \infty |\lambda _j | < \infty$. Note that
the sum converges in $L^1(d\sigma_\epsilon)$. 
The $H^1( d\sigma_\epsilon)$-norm is defined by
%both in the sense
%of distributions (i.e.
%$\int_{\partial\Omega}(\sum_{ j=1}^n \lambda _j a _ j\psi -f\psi) d\sigma
%\, \to 0,$ $\psi \in C^\infty_0(??)$) and in the
% $ H^1 ( d\sigma_ \epsilon)$-norm.  
%The latter is given by
%with respect to the norm 
\begin{equation}\label{EH11}
\| f\|_{H^1(d \sigma_\epsilon)} = \inf \{ \sum|\lambda _j | \}
%\approx \|f^*\|_{L^1(\partial\Omega, d\sigma_{\epsilon})} 
\end{equation}
where the infimum is taken over all possible representations of $f$.
Note that while  $H^1(d\sigma_\epsilon)$-atoms are defined  for
$\epsilon > -1$, 
 we consider (and need for our application in Theorem
\ref{TMPNA})  the space   $H^1( d\sigma_\epsilon)$   only for $\e\le 0$. 
This allows us to avoid having to define spaces of distributions on 
Lipschitz graph domains. (See Coifman and Weiss \cite{CW:1976} or
Stromberg and Torchinsky \cite{ST:1989} for a discussion of these
spaces.)

%For future reference, we also define spaces on subsets of the
%boundary, $N$ and $D$. Here, we are in familiar territory. 
%Chang, Krantz and Stein \cite{CKS:1993}  introduce
%Hardy spaces on domains and use them in 
% other problems in 
% partial differential equations. 
%Also, see the 
%the dissertation of Sykes \cite{JS:1999} 
%(see also Sykes and Brown, \cite{SB:2001}). 
%A function $f$ is in $H^1 (N,  \sigma_\epsilon)$ if and only if $f$ is
%the restriction to $N$ of a function in $H^1 (\partial \Omega,
%\sigma_\epsilon)$. One can  see that such functions can be written as
%sums of atoms where the atoms for $N$ are the restriction to $N$ of
%atoms for $ \partial \Omega$. 
%
%We follow the same path for the Dirichlet data. 
%This data will be taken from the space 
%$H^{1,1}(\sigma_\epsilon)$ for  the regularity problem and from
%$H^{1,1}( D, \sigma_\epsilon)$  in the mixed problem. 

Let $\e\le 0$. We say that $A$ is 
an $H^{1,1}(d\sigma_\epsilon)$-{\em{atom for}} $\partial\Omega$ if for
some $x_0$ in  $\partial \Omega$
$$
A(x)= \int_{x_0}^x  a(t) d\sigma(t), \qquad x \in \Omega 
$$
where $a$ is an $ H^1(d\sigma_\epsilon)$-atom, and $\int_{x_0}^x$ denotes 
 integration along the portion of the 
boundary  $\p\Om$ with endpoints $x_0$ and $x$. 
Given an  $ H^1(d\sigma_\epsilon)$-atom $a$, the integral above defines
$A$ uniquely up to an additive constant.
%\comment{
%; we select a representative
%$A$ 
% which is compactly supported
%in a surface ball $\Delta_r(x)$; for this choice  we have
% $\|dA/d\sigma\|_{\infty}\le 1/\sigma_{\e}(\Delta_r(x))$.
%}
For $ \epsilon \leq 0$, we define the space
$H^{1,1}(d\sigma_\epsilon)$ to be sums of the form  
$$
F(x) = \sum_{ j=1}^ \infty\lambda_j A _j  (x) 
$$
where the coefficients satisfy $ \sum_{j=1}^\infty |\lambda_j| <
\infty$. The norm of $F$ in $H^{ 1,1}(\sigma_{ \epsilon'})$, 
$\|F\|_{H^{1,1}(d\sigma_\epsilon)}$, is defined  to be the infimum of
$  \sum_{j=1}^ {\infty}|\lambda_j|$ over all 
possible representations of $F$ as sums of atoms.  
If we choose the base point $x_0$ for each atom to be
zero, then we have  that the sum defining $F$ converges in the norm given by
$$
\sup |x|^ { -\epsilon } |F(x)| .
$$
Furthermore, it is easy to see that we have 
$$
\|F\|_{ H^ { 1,1}( d\sigma_\epsilon )} = \|\frac { dF}{ d\sigma}\|_{H^
  1(d  \sigma _ \epsilon ) }.
$$
%%46.3
\comment{
The space $H^{1,1}( d\sigma_\epsilon)$, for $-1<\e\le 0$
 is defined 
by taking sums $F=\sum_j \lambda_j a_j$ of
 $H^{1,1}(d\sigma_\epsilon)$-atoms that
converge in the norm 
$$ \sup |x|^{-\e}|F(x)|.$$
Specifically, $F\in H^{1,1}( d\sigma_\epsilon)$ if and only 
if $F=\sum_{j=1}^{\infty} \lambda_j A_j
=\sum_{j=1}^{\infty} \lambda_j \int_{x_0}^x a_j(t)d\sigma(t)$, 
with $\sum_{j=1}^{\infty}|\lambda_j|
<\infty$.S
Since $\sum_j \lambda_j a_j$ converges in $L^1(d\sigma_{\e})$
 the latter yields
$$
F=\int_{x_0}^x \sum_{j=1}^{\infty} \lambda_j a_j(t) d\sigma(t)=\int_{x_0}^x 
f(t) d\sigma(t),
$$
where $f=\sum_{j=1}^{\infty} \lambda_j a_j \in H^1(d\sigma_\epsilon)$. 
The $H^{1,1}( d\sigma_\epsilon)$-norm is defined by
\begin{equation}\label{EH111}
\| F\|_{H^{1,1}(d \sigma_\epsilon)} = \inf \{ \sum|\lambda _j | \}
\end{equation}
where the infimum is taken over all possible representations of $f$.
We have
\begin{proposition}
$$\| F\|_{H^{1,1}(d \sigma_\epsilon)} = \left\| \frac{d F}{d\sigma}\right\|_{H^1
(d \sigma_\epsilon)}.$$
\end{proposition}
The proof is standard and is omitted.
}
%and
%$H^{1,1}(D, \sigma_\epsilon)$ in the same way that constructed $H^1$.
By the {\em Neumann problem for} $H^1(d\sigma_\epsilon )$ we mean the
Neumann problem for $L^1( d\sigma_\epsilon)$, see
 (\ref{Neumann}), where the data $f_N$ is now taken from
 $H^1( d\sigma_\epsilon)$.
% (and the boundary data convergence is
% interpreted in the sense of distributions). 
By the {\em regularity problem for} $H^{1, 1}(d\sigma_\epsilon )$ we mean the
regularity problem for $L^1( d\sigma_\epsilon)$, see
 (\ref{Regular}), where the data $f_D$ is taken from
 $H^{1,1} ( d\sigma_\epsilon)$. 
%(and the boundary data convergence
% is again 
% to be interpreted in the sense of distributions).

For future reference, we also define spaces on subsets of the
boundary. 
%Here, we are in familiar territory.
%Chang, Krantz and Stein \cite{CKS:1993}  introduce
%Hardy spaces on domains and use them in
% other problems in
% partial differential equations.
%Also, see the
%the dissertation of Sykes \cite{JS:1999}
%(see also Sykes and Brown, \cite{SB:2001}).
A function $f$ is in $H^1 (N,  d\sigma_\epsilon)$ if and only if $f$ is
the restriction to $N$ of a function in $H^1 (d\sigma_\epsilon)$. Such functions can be written as
sums of atoms that are restrictions to $N$ of $H^1(d\sigma_\epsilon)$-
atoms for $ \partial \Omega$. The space $H^{1,1}(D,d\sigma_\epsilon)$
is defined in a similar fashion. See Chang, Krantz and Stein \cite{CKS:1993},
 Sykes \cite{JS:1999}, and Sykes and Brown \cite{SB:2001} for 
 additional works that study Hardy spaces on domains. 
 
%
%We follow the same path for the Dirichlet data.
%This data will be taken from the space
%$H^{1,1}(\sigma_\epsilon)$ for  the regularity problem and from
%$H^{1,1}( D, \sigma_\epsilon)$  in t

The main result of this section is the following theorem. We omit the proof
since it is quite similar to the argument for the mixed problem in the
following section (Theorem \ref{TMPNA}). 

\begin{theorem} \label{rnh1} Let $\Omega$ be a standard
 Lipschitz graph domain with Lipschitz constant $M$ and
 let $\e_0(M)<\e'\le 0$ be small.
%|\epsilon'| \leq \epsilon_0(M)$ be small. 
Then the $H^1 (d\sigma_
   {\epsilon'})$-Neumann and the
 $H^{1,1}(d\sigma_{\epsilon'})$-regularity
problems are uniquely solvable. The solutions satisfy, 
 respectively,
\begin{eqnarray*}
\|  u  \| _{ H^{1,1} ( d\sigma_{\epsilon'})}+
\int_{\partial \Omega} ( \nabla u )^* (x) \, d\sigma_{\epsilon'} & \leq & C
\left\| \frac {\partial u }{ \partial \nu}\right\| _{ H^1 (d\sigma_{\epsilon'}) }
\\
\left\| \frac {\partial u }{ \partial \nu}\right\| _{ H^1
  (d\sigma_{\epsilon'}) }+ 
\int_{\partial \Omega} ( \nabla u )^* (x) \, d\sigma_{\epsilon'} & \leq
&  C
\|  u  \| _{ H^{1,1} ( d\sigma_{\epsilon'})}.
\end{eqnarray*}
\end{theorem}

\begin{remark} 
We do not assert uniqueness when $ \epsilon' >0$. 
In the rest of this  paper, we will only use the existence of a
solution for the  $H^1(d\sigma_{\epsilon'})$-Neumann problem in the case
when $\epsilon'\leq 0$. 
(In order to fully treat the case $\epsilon'>0$, one needs to give a different
definition of the normal derivative at the boundary.
For  $\epsilon'>0$ the $H^1(d\sigma_{\epsilon'})$-boundary data may
not be in $L^1_{loc}(d\sigma)$ and hence fail to be a function.   
%In
%the rest of this paper, we will only use the existence when $\epsilon
%>0$. Thus, we will not establish uniqueness when $ \epsilon >0$.
See Brown \cite{RB:1995a} and Fabes, Mendez and Mitrea \cite{FMM:1998}
for a treatment of the Neumann problem with data which is not 
locally integrable).
%functions. 
%As these works indicate, we need to make a new definition
%of $\partial u/\partial \nu$ at the boundary for such functions.
\end{remark}

\section{The mixed problem with atomic data. }
\setcounter{equation}{0}

By the \textit{mixed problem for } $H^1(d\sigma_{\epsilon'})$ we mean the 
mixed problem for $L^1( d\sigma_{\epsilon'})$, see
(\ref{Mixed}), where the Dirichlet data $f_D$ is taken from
$H^{1,1}(D, d\sigma_{\epsilon'})$ and the Neumann data $f_N$ is in
$H^1(N, d\sigma_{\epsilon'})$.
In this section, we consider the mixed problem where the Neumann data
is an $H^1(N,  d\sigma_{ \epsilon '})$ atom. 
and the Dirichlet data is zero.
Since atoms lie in $L^2(d\sigma_\epsilon)$ for $ \epsilon > -1$, by
Theorem 4.6 we have existence and uniqueness of the solution to the
mixed problem for $L^2(d\sigma_\epsilon)$ with these data. 
Our first goal is to show that this solution also has non-tangential
maximal function in $L^1(d\sigma_{ \epsilon'})$ for $\epsilon'$ near
zero. 
%%
%% Since atoms lie in $L^2(d
%%\sigma_\epsilon)$ for $ \epsilon > -1$, we  have the existence of the
%%$L^2(d\sigma_\epsilon)$  solution with this data.
%%Our goal in this section is to show that the non-tangential
%% maximal function of the gradient of the
%% $L^2(d\sigma_\epsilon)$-solution also lies in
%% $L^1(d\sigma_{ \epsilon '})$, see Theorem \ref{53.5}.  
%%
%%This is the key step in proving the existence of solutions with data
%%in $H^1$ and also in $L^p$.  
% However,
%considering results in 
% sectors indicates that we do not expect to have results in
%$L^2(\sigma)$, even in smooth domains.  \marginpar{Need ref?}
%Our first  goal is to prove the following theorem.
%We will then interpolate between this result and the result in $L^2(
%\sigma_\epsilon)$ to obtain estimates in $L^p(\sigma)$. 

\begin{theorem} \label{53.5}  Suppose $\Omega$, $N$, $D$ is a standard graph
domain for the mixed problem with Lipschitz constant $M<1$ and set $
\beta = \tan^{ -1} M$.
% satisfying the conditions of Lemma%\ref{LVFMP},
Then, there is $\delta= \delta (M)$ with $ 1> \delta>0$ so that, for 
%$2\beta/(\pi - 2\beta)<\epsilon<1$
% (where we have set $\beta = \tan^{-1}M>0$) and
$\epsilon '$ satisfying  $\frac{4\beta-\pi}{2(\pi-\beta)}<\epsilon'<\delta$
%(where we have set $\beta = \tan^{-1}M>0$),
 we may solve the mixed problem
% there is an $\epsilon_0 = \epsilon_0 (M)$ so that, if $ |\epsilon | <
%\epsilon_0$ and $a$ is an $H^1 (N, d\sigma_{ \epsilon'})$-atom,
%we may solve
(\ref{Mixed}) for $H^1( d\sigma_{\epsilon'})$ with zero
Dirichlet data and with Neumann data an
$H^1 (N, d\sigma_{ \epsilon'})$-atom, $a$.
% $a$
% and with zero
%Dirichlet data.
 The solution $u$ satisfies the 
estimate
\begin{equation}\label{EMPH11}
\int_{ \partial \Omega} (\nabla u )^* (x) \, d \sigma_{ \epsilon'}
\leq C(M, \epsilon')\, .
\end{equation}
In addition, we have
\begin{equation}\label{EMPH12}
   \|u \| _{ H^{1,1}( d\sigma_{\epsilon'})} +
   \left\|\frac {\partial u } {\partial \nu}\right\|_{ H^1
   (d\sigma_{\epsilon '}) } \leq  C(M, \epsilon'). 
\end{equation}
\end{theorem}

As a step towards the proof of Theorem \ref{53.5}, we construct a
Green's function for the mixed problem using the
method of reflections--an old idea that was used by Dahlberg and Kenig
\cite{DK:1987} 
to obtain a similar result for the Neumann and regularity problems.
The estimates for the Green's function are a consequence of the
H\"older regularity of weak solutions of divergence-form
equations with bounded measurable coefficients.

We begin by constructing a bi-Lipschitz map, $\Phi:\reals^2
\rightarrow \reals ^2$ with $ \Phi(\Omega) =Q$ where $Q$ is the first
quadrant, $Q=\{(x_1,x_2) : x_1>0, x_2>0\}$. We also require that 
$
\Phi(N) = \{ (x_1, x_2) : x_1 \geq 0, x_2=0\}$ and 
$\Phi(D)= \{(x_1,x_2) : x_1=0, x_2>0\}$.
On $Q$ 
we now define the operator $L=\mathrm{div}A\nabla$ where the coefficient
 matrix $A$ has bounded entries (these are first-order derivatives of $\Phi$),
so that  $u$ satisfies
$$
\left \{  \begin{array}{ll} \Delta u = 0, \qquad  &\mbox{in }\Omega\\
\frac {\partial u }{ \partial \nu } = f_N, \qquad  & \mbox{on } N \\
u=f_D, \qquad & \mbox{on } D \\
	  \end{array}
\right.
$$
if and only if the function $v$ defined by  $v=  u\circ \Phi^{-1}$
  satisfies
$$
\left\{\begin{array}{ll}
Lv= 0 , \qquad & \mbox{in } Q\\
A \nabla v \cdot \nu = f_N\circ \Phi^{ -1} & \mbox{on } \Phi(N)\\
v= f_D \circ \Phi^{ -1} \qquad & \mbox{on } \Phi(D)
       \end{array}
\right. 
$$
Now, we extend the coefficients of $L$ by reflection,
 so that $Lv=0$ if and only if
$ L(v\circ R_j)=0$ where $R_1$ and $R_2$ are the reflections
$
R_1(x_1,x_2)= (-x_1, x_2)$ and  $R_2(x_1,x_2) = ( x_1, -x_2)$.
Next, letting $G$ denote the Green's function for $L$ in $Q$, we set 
$$\mathcal M(z,w) = G(z,w)-G(z,R_1w)+ G(z,R_2w) -G(z,R_1R_2w)\, ,\quad z, w\in Q
$$ and observe that $\mathcal M$ is a Green's function for the mixed
problem in $Q$. 
Recall that  we may
find a Green's function for $L$ in $Q$ which satisfies
$$
|G(z,w) | \leq C ( 1+ |\log |z-w||)\, ,\quad z, w\in Q
,
$$
see Kenig and Ni \cite{KN:1985}.
We then observe that if $|z-\zeta | =1$ and $ |\zeta- w| < 1/2$, then 
$|G(z,\zeta) -G(z,w)| \leq C | \zeta-w|^ \delta$, where
$0<\delta<1$, $\delta=\delta(M)$. This is a standard
estimate of H\"older continuity for solutions of divergence-form
elliptic operators.  
Finally, by rescaling, we obtain  
\begin{equation} \label{GreenHolder} \label{54.3}
|G(z,\zeta) -G(z,w)| \leq C \left( \frac{ | \zeta- w|
  }{|z-\zeta|}\right)^\delta , \qquad \mbox{if } |\zeta- w| < \frac 1 2 |
z-\zeta|.
\end{equation}
The latter immediately implies the same estimate for $\mathcal M$ in $
Q$, namely 
%%55.5
\begin{equation}\label{MHolder}
|\mathcal M(z,\zeta) -\mathcal M(z,w) | \leq C \left(\frac  {|\zeta- w| } {
 |z-\zeta|}\right)^\delta \quad \mbox{if } z , 
\zeta,w \in \overline{Q}\  
\mbox{ and }|\zeta-w|< \frac 1 2 | z-\zeta|
\end{equation}

We will need an additional estimate for $\mathcal M(z,\zeta)$ when
 $\zeta \in Q$ is near
 $\Phi(D)$. Since we have
$\mathcal M(z,x) = 0$ if $z\in Q$ and $x \in \Phi(D)$, it follows by
% the 
% continuity of $M$ 
continuity that ${\cal M}$ is small near $\Phi(D)$. More precisely, 
let $\zeta \in Q$ and 
suppose that 
$\hat x$ is a point on $\Phi(D)$ for which $|\hat x -\zeta | =
 \dist ( \zeta , \Phi(D))$.
 For $z\in Q$, we have $M (z,\hat x)= 0$ and 
%the H\"older continuity of $M$ in
(\ref{MHolder}) implies 
\begin{equation}\label{MDHolder}
| \mathcal M(z,\zeta) |\leq C \frac {|\zeta|^ \delta }{|z-\zeta|^ \delta }\, , 
\quad \mathrm{if}\ z, \zeta \in Q\ \, \mathrm{and}\ 
\mbox{dist}(\zeta ,D) <\frac 1 2 |z-\zeta |,
\end{equation}
(here we have used the fact that $0\in \partial Q$).  
By using the change of variables $\Phi$, it is easy to translate back
to the original domain $\Omega$. Since $ \Phi$ is bi-Lipschitz we
obtain estimates (\ref{MHolder}) and (\ref{MDHolder}) in $ \Omega$. 

%These estimates will give the behavior at infinity of solutions to the
%mixed problem which have atomic Neumann data.
 These estimates are a key
ingredient in the study of
%the next Lemma which describes
 the behavior
of the solution of the mixed problem with atomic
Neumann data. More precisely, we have

\begin{lemma}\label{58.5} Assume $\epsilon'>-1$. 
Let $u$ be a solution of the mixed problem (\ref{Mixed}) for
 $L^1( d\sigma_{\epsilon'})$, where the
 Neumann data
is an $H^1(N, d\sigma_{\epsilon'})$-atom, $a$,
 and the Dirichlet data is zero.
 Then, for any integer $k\ge 1$,
$u$ satisfies   
$$
|u(z) | \leq C \frac { \rho^\delta}{ |z-x_a|^ \delta} \frac {
  \sigma ( \Delta _ \rho ( x_a) ) }{ \sigma_{ \epsilon'} ( \Delta _
  \rho ( x_a))}, \quad z\in \Omega\, ,\quad |z-x_a|\geq 2^k \rho.
$$
Here, $\delta > 0$ is as in the estimate for the Green's function for
the mixed problem (\ref{54.3}) and $\Delta_\rho(x_a)$ is the surface 
ball where
$a$ is supported.
\end{lemma}

\begin{proof} We consider two cases:
1) $\Delta _ \rho ( x_a) \subset N$ and 2) $\Delta _ \rho ( x_a) \cap
   D \neq \emptyset$. 

In case 1), we have 
$$
u (z) = \int_N \mathcal M(z,x) a(x) \, d\sigma(x) 
=\int_N (\mathcal M(z,x)-\mathcal M(z,x_a)) a(x) \, d\sigma(x)\, ,\quad z\in \Omega 
$$ 
where the second identity uses the fact that the atom $a$ has mean value zero.
Next, we use the continuity of ${\mathcal M}$ from (\ref{MHolder}) to
obtain 
\begin{equation}\label{58.6}
|u(z) | \leq C \frac { \rho ^ \delta}{|z-x_a|^ \delta }\int _N |a|\,
 d\sigma, \quad \mathrm{for}\ z\in \Omega\, ,\ |z-x_a| > 2^k \rho .
\end{equation}
Finally, the normalization of $a$ in the definition of an atom implies
that 
$$
\int_N |a| \, d\sigma \leq \frac { \sigma  ( \Delta _ {\rho}  (x_a)  )}
{\sigma_{ \epsilon' }( \Delta _\rho ( x_a))}\, .
$$
%provide $ \epsilon' > -1$. 
This completes the proof in case 1.

In case 2), we do not have $\int_N a \,d\sigma =0$. However, 
 estimate (\ref{MDHolder}) yields 
$$
|u(z) | \leq C \frac { \rho ^ \delta}{ |z-x_a|^ \delta} \int _N |a| \,
d\sigma \, ,
$$
and then  we use the normalization of $a$ to conclude the proof. 
\end{proof}

Before proceeding, we need a few technical results.  In this Lemma and
below, given a point $x_a \in\partial\Omega$, we consider a 
 ball $B_{\rho}(x_a)$ and a boundary ball
 $ \Delta _ \rho ( x_a)$, and set:
%then construct
 $ R_k = B_{ 2^ {k+1}\rho }( x_a)\setminus
B_{2^k \rho }( x_a)$, $ \tilde R_k = B_{ 2^{k+2} \rho}( x_a) \setminus
B_{ 2^{ k-1} \rho } ( x_a)$, $ \Delta _k = \sball {x_a}{ 2^k \rho }$
and  $ \Sigma _k = \Delta _{k+1}\setminus \Delta _k $. 
With these definitions, we can now  state 
\begin{lemma} \label{61}
 Let $\lambda\in\mathbf R$
 and set $\alpha(z) = (\Re ( e^{i\lambda} (z)^ \epsilon), 
 \Im ( e^{i\lambda} (z)
  ^ \epsilon))$. Then, we have
$$
\int_{R_k }|\alpha | ^ p \, dy 
\leq C 2^ k \rho \int _{ \Sigma_k } | \alpha |^ p \, d\sigma 
\leq C 2^ k \rho \sigma_{ \epsilon p }( \Sigma _k )
$$
provided $ \epsilon p > -1$.
\end{lemma}

\begin{proof} The proof is a computation and  we omit the details.
\end{proof}
We now take a brief detour to discuss solutions of the mixed problem
in the energy sense. 
%We suppose that $ B$ is a ball with center in $ \Omega$ and set $A= B
%\cap \Omega$. 
Let $B$ be a ball with center in $ \Omega$. 
We say that $u$ is an  {\em energy solution of the mixed problem in
 }$B$, 
$$
 \left\{
\begin{array}{ll}
\Delta u = 0  \qquad & \mbox{in } B\cap \Omega \\
u = 0  \qquad & \mbox{on } D\cap B  \\
\frac { \partial u } {\partial \nu } = 0  \qquad & \mbox{on } N \cap B\,  
 %\right. 
\end{array}\right.
$$
%$u$ lies in the Sobolev space $W^{1,2} (B\cap\Omega)$
if  $u$ lies in the Sobolev space
$W^{1,2} (B\cap\Omega)$, $u$ vanishes on $D\cap B$  and for
every $v$ that lies in $W^{ 1,2}(B\cap\Omega)$ and vanishes on $
\partial(B\cap\Omega)  \setminus N$ we have
$$
\int _{B\cap\Omega} \nabla u \cdot \nabla v \, dy = 0.
$$
Using a Carleson measure argument, see Proposition \ref{Carleson},
it is not difficult to see that if $u$ is a solution of the mixed problem
in $L^2(d\sigma_{\e})$ then $|\nabla u|^2$ is integrable
on  bounded subsets of $\Om$, provided $\e<1$. This shows
that  a solution for the mixed problem in $L^2(d\sigma_{\epsilon})$
is, in particular, an energy solution.

In the next lemma, we use  $-\!\!\!\!\!\int_E f \, dx$ to denote the
average, 
$-\!\!\!\!\!\int_E f(x) \, dx \assign |E|^{-1} \int f(x) \, dx$. 
%%\marginpar{Is there a restriction on $r$.}
\begin{lemma}  \label{61.5}
Let $\Omega$, $D$, $N$ be a standard Lipschitz graph domain for the mixed
 problem with Lipschitz constant $M$. There is an exponent
$q_0=q_0(M)>2$  so that on any  ball $B$ with center in $ \Omega$, if $u$ is
an energy solution of
$$\left\{
\begin{array}{ll}
\Delta u = 0, \qquad & \mbox{in } 2 B \cap \Omega \\
u = 0, \qquad & \mbox{in } D \cap 2B\\
\frac { \partial u } { \partial \nu } =0 , \qquad &  \mbox{on } N \cap 2B
\\
\end{array} \right. 
$$
then, for all $1\le q\le q_0$,  we have 
$$
\left( \average _{ B\cap \Omega } |\nabla u |^ q \, dx \right)^ { 1/q} 
\leq \frac C r \left( \average _{ 2B \cap \Omega } |u |^2 \, dx \right
) ^ {1/2} ,
$$
 where $r$ is the radius of $B$.
\end{lemma}

\begin{proof}
This follows from the Caccioppoli inequality and a reverse H\"older
inequality as in  Giaquinta \cite{MG:1983}.
\end{proof}
%61.51

The next estimate gives the decay at infinity of a solution to the
mixed problem with atomic data. 

\begin{lemma} \label{62}
Let $\Omega$, $D$, $N$ be a standard Lipschitz graph domain for the mixed
 problem (\ref{Mixed}).
 Let $\epsilon'>-1$ and $-1+2/q_0 <\epsilon<1$, where $q_0$
 is as in Lemma \ref{61.5}. 
If $u$ is a solution of the mixed problem for  $L^2 (
d\sigma_ \epsilon)$  with zero Dirichlet data and
with Neumann data  a $ \sigma_{ \epsilon'}$-atom, $a$, which is
supported in $ \Delta _\rho ( x_a)$, then we have
$$
\int _{ \Sigma_k } |\nabla u | ^2 \, d\sigma_ \epsilon
\leq C \frac {\sigma_\epsilon ( \Sigma_k )}  {( 2^k \rho ) ^2 } \left[
  \sup _{ \cup _{|j|\leq 1} 
 \tilde R_{ k+j}} |u| \right]^2\, ,\qquad \mathrm{for \ all}\ k\geq 2\, . 
$$
%Here, we assume the Lipschitz constant $M<1$ and $ \epsilon $ is as in
% Lemma \ref{LVFMP}.
\end{lemma}

\begin{proof} 
%We let $ \alpha$ be a vector field 
%%with $ |\alpha (x) |
%%=  |x|^ \epsilon$ as constructed in Lemma \ref{LVFMP}. 
%of the form $\alpha(z) = e^{i\lambda}z^{\epsilon}$, for some
%$\lambda \in \mathbf R$ and for $\epsilon$ as in the hypothesis.  
With the same notations as in Lemma \ref{61}, we let $
  \eta _k\geq 0$ be a smooth cutoff function which equals one on $R_k$,
is   zero outside $ \tilde R _k$, and satisfies $|\nabla \eta|\leq
  C/(2^k\rho)$. 
We let $ \alpha$ denote a vector field
%with $ |\alpha (x) |
%=  |x|^ \epsilon$ as constructed in Lemma \ref{LVFMP}.
of the form $\alpha(z) = (\Re ( e^{i\lambda}z^{\epsilon}), \Im (e^{ i
  \lambda} z^ \epsilon))$, for some
$\lambda \in \reals $ and for $\epsilon$ as in the hypothesis.

 We apply the Rellich identity with
  vector field $ \alpha \eta_k$ as in Lemma \ref{RellichIdentity} to conclude
  that for $k \geq 2$, we have 
\begin{equation}\label{62.1}
\int _{ \Sigma_k } |\nabla u | ^ 2 \, d \sigma _ \epsilon
\leq 
\frac C { 2^k \rho} \int _{ \tilde R_k } | \nabla u | ^2 |\alpha | dy.
\end{equation}
(This uses the fact that the mixed data of $u$ is zero on the support
 of $ \eta
_k$ when $k \geq 2$).  

Applying H\"older's inequality 
%for $p \geq 1$ ($p$ to be determined later),
we obtain
\begin{equation} \label{62.2} 
\int_{ \tilde R_k } |\nabla u |^2 |\alpha | \, dy
\leq \left( \int_{ \tilde R_k } | \nabla u | ^{ 2p} \, d x\right ) ^
     {1/p}
\left ( \int _{ \tilde R _k} |\alpha |^{ p'} \, dy \right) ^ { 1/p'},
\end{equation}
where $1/p + 1/p' =1$ and $p$ is lies in the interval
$1/(1-|\epsilon|)<p\le q_0/2$ (if $-1+2/q_0<\epsilon<0$) or, if $\epsilon>0$,
$1\le p\le q_0/2$. (These conditions grant that Lemmata \ref{61.5}, 
 \ref{61} and \ref{measure-property} apply in what follows).
We now cover $\tilde R_k$ with a (fixed) number of balls $B_{k,n}$ (each
 centered at a point in $\Omega$), $n=1,...,m$, such that
 $\mathrm{diam} B_{k,n} \approx 2^k\rho$ and
$$
\tilde R_k \subseteq \cup_{n=1}^m B_{k,n} \subseteq \cup_{|j|\le 2} R_{k+j}
 = \cup_{|j|\le 1}\tilde R_{k+j}\, .
$$
By Lemma \ref{61.5}, for $p$ as above, we have
$$
\left( \int_{ \tilde R _k } |\nabla u |^ {2p} \, dy \right) ^ {1/p}
\leq C
( 2^ k \rho) ^ { \frac 2 p - 4} \sum _{|j | \leq 1 }  \int _{\tilde R_{k+j}}
  |u|^ 2 \, dy. 
$$
Moreover, Lemma \ref{61} and Lemma \ref{measure-property}, also for
 $p$ as above, imply
\begin{equation}\label{62.3}
\left( \int _{\tilde R_k} |\alpha |^{p'} \, dy \right) ^{1/p'} 
\leq C ( 2^k \rho ) ^ { 1-\frac 2p} \sigma_\epsilon ( \Sigma_k ).
\end{equation}
%provided $ \epsilon p'>-1$.
Combining (\ref{62.1}) to (\ref{62.3}), we obtain 
\begin{eqnarray*}
\int_{\Sigma_k } |\nabla u | ^2 \, d \sigma _ \epsilon
%& \leq & \frac C { 2^ k \rho} \int _{ \tilde R _k } | \nabla u |^ 2
%|\alpha | \, dx \\ 
%%& \leq  & 
%%\frac C { 2^ k \rho }  \int _{ \tilde R_k } |\nabla u | ^2 |\alpha |
%%\, dx  \\
%&\leq &
% \frac C { 2^ k \rho} ( 2^ k \rho) ^ { \frac 2 p - 4} 
%\sum _{ |j|\leq 1 } \int _{ R_{ k+j} } | u | ^ 2 \, dx \left(
%  \int_{R_{k+j}} |\alpha |^ {p'} \, dx \right ) ^{ 1/p'} \\
& \leq &
 \frac C { 2^ k \rho }( 2^ k \rho ) ^ { \frac 2 p - 4 } ( 2^ k
\rho ) ^ { 1- \frac 2 p } \sigma _ \epsilon(\Sigma_k) 
%( \Delta _{ 2^ k \rho })
\sum _{ |j | \leq 1 }  \int _{\tilde R_{ k+j }}|u |^ 2 \, d x \\
& \leq & 
C \sigma_ \epsilon ( \Sigma _k ) ( 2^ k \rho ) ^ { -2}
\left[ \sup _{ x \in \cup _{ |j | \leq 1 } \tilde R _{ k + j } } | u |
  \right]^ 2 \, .
\end{eqnarray*}
\end{proof}
%62.5

\begin{lemma}\label{63}  Let $\Omega, N, D$ be a standard Lipschitz graph
 domain for the mixed problem. Suppose $u$ is the solution of the mixed
 problem
 for $L^2( d\sigma_\epsilon)$, see (\ref{Mixed}), with zero Dirichlet data and
  with Neumann data an
  $H^1(N, d\sigma_{\epsilon '})$-atom, $a$, supported in a surface ball
$\Delta_{\rho}(x_a)$. Let $\epsilon'>-1$ and $-1+2/q_0<\epsilon<1$
 (with $q_0$ as in Lemma \ref{61.5}). Then, for all integers $k\ge 5$
 we have
 % Suppose Lemma \ref{62} is valid for $ \epsilon$ and Lemma \ref{58.5}
 % is valid for $ \epsilon'$. Then we have that
$$
\int_{ \Sigma_k } ( \nabla u )^ *(x) ^2 \, d\sigma_ \epsilon 
\leq 
\frac { C \sigma_ \epsilon ( \Sigma_k )}
{2^{2k(1+\delta)} ( \sigma_{ \epsilon '}( \Delta_\rho ( x_a)))^2},
$$
where $ \delta $ is as in (\ref{GreenHolder}).
\end{lemma}

\begin{proof} By Lemmata \ref{58.5} and \ref{62}, we have 
\begin{eqnarray}
\int _{ \Sigma _k } |\nabla u |^2 \, d\sigma _ \epsilon 
 & \leq &  \frac C { (2^k \rho )^2} \sigma_ \epsilon ( \Sigma _k )
 \left(\frac {\rho ^ \delta}{ ( 2^k \rho )^ \delta } 
\frac { \sigma ( \Delta _ \rho ( x _a))}{\sigma_{\epsilon '} ( \Delta _
  \rho ( x_a))}\right)^2 \nonumber \\
 & = & \frac C {(2^k \rho)^2}\sigma_\epsilon ( \Sigma _k )2 ^ { -2k\delta}
\left(\frac  \rho  {\sigma_{\epsilon '} (\Delta_\rho (
  x_a))}\right)^2. \label{63.1}
\end{eqnarray}

In order to pass from the estimate above for $ \nabla u$ to an estimate for
$(\nabla u)^*$ , we use the Cauchy
kernel to represent $ \partial u$, the complex derivative of $u$, see
 (\ref{ECderiv}). 
%of course, a similar argument
%works for $ \bar \partial u$. 
To carry out this argument, we 
%fix $x \in \Sigma_k$ and 
%let 
consider a cutoff function $ \eta_k$
% be %a cutoff function which 
that is supported in $B_{ 2^{ k+4}\rho}(x_a)
    \setminus B_{ 2^{k-3}\rho }(x_a)$,
equals one on  $  B_{ 2^{k+3} \rho }(x_a) \setminus
  B_{ 2^{ k-2}\rho }(x_a)$ and, furthermore, satisfies:
$|\partial\eta_k|\le C/2^k\rho$.
 On account of the analyticity of $\partial u$,
the Cauchy-Pompeiu formula, see e.g. Bell \cite{SB:1992}, 
 yields 
\begin{equation} \label{ShenRep}
\eta_k(z)  \partial u (z) =  K_{\partial\Omega}(\eta_k\partial u)(z) + I(z),
\quad z\in \Omega, 
\end{equation}
where we have set 
%\frac 1 { \pi  } \int _{ \Omega } 
%\frac 1 { z-y } \bar \partial \eta^x (y) \partial u(y) \, dy
\begin{equation}\label{ShenRep1}
K_{\partial\Omega}(\eta_k\partial u)(z)\, =\, 
 \frac 1 { 2\pi i  } \int _{ \partial \Omega } \frac 1 { z-\zeta}
 \eta_k (\zeta)
\partial u ( \zeta )\, d\zeta\, ,\quad z\in \Omega, 
% ( \nu_1(y) + i \nu _2(y)) \, d \sigma (y ). 
\end{equation}
and   $d\zeta$ is used to  denote complex line integration. 
Also, we define 
\begin{equation} \label{ShenRep2}
I(z) =
\comment{
\frac 1 { \pi  } \int _{ \Omega }
\frac 1 { z-y } \bar \partial \eta_k (y) \partial u(y) \, dy\,
}
\, \frac 1 { \pi  } \int _{R_{k-2}\cup R_{k+3}}
\frac 1 { z-y } \bar \partial \eta_k (y) \partial u(y) \, dy\, 
\end{equation}
and in this expression, $dy$ denotes area measure. 

Next, for $x\in\Sigma_k$
we decompose the sector  $\Gamma (x) 
= \Gamma ^n (x) \cup \Gamma^f (x)$ where
$\Gamma ^n (x) = \Gamma (x) \cap B_{ \kappa 2 ^k \rho }(x_a)$,
$ \Gamma^ f (x) =\Gamma(x) \setminus B_{ \kappa 2^k \rho }(x_a)$, and 
the constant $ \kappa $ is chosen so that
% if $ x $ is in $
%\Sigma _k $, then 
$ \Gamma ^n (x) \subset \cup _{ |j|\leq 1} R_{
  k+j}$ for  $x \in \Sigma_k$. 
If we let $ v_n^*(x)$ denote the supremum of $|v|$ on $ \Gamma _n (x)$ and
similarly for $v_f^*$, using (\ref{ShenRep}) and
 the theorem of Coifman, McIntosh and
Meyer \cite{CMM:1982}  we obtain 
\begin{eqnarray} \label{NearEstimate}
\int _{ \Sigma _k } (\nabla u)_n^*(x)^2 \, d \sigma _\epsilon(x)
 & =&   
\int _{ \Sigma _k }  ( \partial u )_n ^ * (x)^2 \, d \sigma _\epsilon
(x) \\ \nonumber
& \leq& 
\sum _{ |j|\leq 3 } \int _{ \Sigma _{ k+j} } | \partial  u (x) |
^ 2 \, d\sigma_\epsilon (x) \\ \nonumber
& & 
\qquad +
\sigma_\epsilon ( \Sigma_k )\left[
\sup_{x\in\Sigma_k}
\left(\sup_{z\in\cup_{|j|\le 1}R_{k+j}}|I(z)|\right)^2
\right] ,
% 
%+ \sigma _\epsilon ( \Sigma _k ) \left
% [ \sup_{ x \in \cup_{|j| \leq 1} R_{ k+j}}
%  I(x) \right ]^2 
\end{eqnarray}
%Here, $ I^x(z)$ is the first integral on the right-hand side of
%(\ref{ShenRep}). 
where we have used that $\supp \eta_k \subset \cup_{|j|\leq
 3}R_{k+j}.$ 

% for all $x\in\Sigma_k$.
We now estimate the term $I(z)$: by the Cauchy-Schwartz
 inequality we have 
$$ |I(z)|\,\le\,
\frac{C}{\rho 2^{k-1}}\bigg(\int_{R_{k-3}\cup R_{k+3}}
|\nabla u|^2\, dy\bigg)^{\frac{1}{2}}\, ,\quad
z\in\cup_{|j|\le 1}R_{k+j}\, . 
$$
On account of the vanishing boundary conditions on $u$ and
 $\partial u/\partial\nu$ we may now apply Caccioppoli inequality
and conclude 
$$
I(z)  \leq \frac  1 { 2^ k \rho  }
\sup_{ y \in \cup_{|j| \leq 4} R_{ k+j}} |u(y) |\, \quad \mathrm{for\ all}
\ z\in \cup_{|j| \leq 1} R_{ k+j}\quad \mathrm{and}\  x\in \Sigma_k\, .
$$
Using interior estimates, we also obtain  
$$
( \nabla u )_f^* (x) \leq \frac 1 { 2^ k \rho } \sup _{z\in\Omega
\setminus B_{2^{ k-1}\rho}} 
|u(z) |\, ,\quad x\in\Sigma_k\, . 
$$
  The conclusion now follows from (\ref{63.1}),  (\ref{NearEstimate})
and Lemma \ref{58.5}.  
\end{proof}

We are ready to prove Theorem \ref{53.5}. 
\begin{proof}[Proof of Theorem \ref{53.5}]
We fix $\e>0$ such that $2\beta/(\pi-2\beta)<\e<1$ 
so that, for $\e'$ as in the hypotheses we have: 
$(\e-1)/2 <\e'$.
% As mentioned before, 
%it is enough to consider the case when the
%Dirichlet data is zero. 
%Moreover, we may assume that the Neumann data is a
% $\sigma_{\epsilon'}$-atom $a$
%supported in 
%$\Delta _\rho (x_a)$, $x_a \in \partial\Omega$.
The proof is based on the following elementary observations:
 given any $\epsilon'>-1$ and $\epsilon>-1$, if $a$ is a
$\sigma_{\epsilon'}$-atom then $a$ lies in
$L^2 (N, d\sigma_\epsilon)$; by Theorem \ref{MixedProb2}  
 we may solve the mixed problem (\ref{Mixed}) for $L^2( d\sigma_{\epsilon})$
 with Neumann data $a$ and zero Dirichlet data
 (provided $\beta$ and $\epsilon$ are as in the hypothesis). We 
 let $u$ denote such a  solution and show that $u$ satisfies (\ref{EMPH11})
 and (\ref{EMPH12}). 
We  begin by studying $ \nabla u $ near the support of $a$.
%$ \Delta_\rho (x_a)$. 
By the
Cauchy-Schwarz inequality, together with the  normalization of $a$ and the
estimate for the $L^2 (d\sigma_\epsilon)$-mixed 
problem
%(see Theorem \ref{39.3}),
(see Theorem \ref{MixedProb2}), 
% and the normalization of $a$
we have
\begin{eqnarray*}\label{Elong}
\int _{\Delta_{2^{10}\rho }( x_a)} ( \nabla u )^ * (x) \,
d\sigma_{\epsilon '} 
& \leq &  \left( \int_{ \Delta _{2^{10} \rho } ( x _a)} ( \nabla u )^ * (x)
  ^2 \, d\sigma _\epsilon\right)^ {1/2} \sigma_{2\epsilon'-\epsilon }
  ( \Delta _{ {2^{10} \rho} }( x_a) )^{1/2}  \\
& \leq &  C \| a\| _{ L^2 ( \sigma_\epsilon )} \sigma_{ 2\epsilon '-
  \epsilon }( \Delta _{2^{10} \rho} (x_a)) ^{ 1/2} 
 \\
& \leq &  C \frac { \sigma_ \epsilon ( \Delta _\rho ( x _a))^ {1/2}
  \sigma_{ 2\epsilon'- \epsilon } ( \Delta _{ 2^ {10} \rho} ( x_a) ) ^
	{ 1/2} }
{ \sigma _{ \epsilon ' } ( \Delta _ \rho ( x_a))}\, ,
%\leq C.
\end{eqnarray*}
%65
provided $2\epsilon'- \epsilon >-1$, that is $\epsilon'>(\epsilon -1)/2$
(note that the latter is bounded below by $(4\beta - \pi)/2(\pi-2\beta)>-1$).
Lemma \ref{51} now grants 
%$ \epsilon > -1$, $\epsilon ' > -1$ and $ 2 \epsilon ' -\epsilon > -1$. 
$$
\frac { \sigma_ \epsilon ( \Delta _\rho ( x _a))^ {1/2}
  \sigma_{ 2\epsilon'- \epsilon } ( \Delta _{ 2^ {10} \rho} ( x_a) ) ^
        { 1/2} }
{ \sigma _{ \epsilon ' } ( \Delta _ \rho ( x_a))}\, 
\leq C.
$$
Next, we consider $\nabla u $ away from the support of $a$: 
 we will show that there is $\eta>0$ so that
\begin{equation} \label{64}
\int_{ \Sigma_k } ( \nabla u )^* (x) \, d \sigma_{ \epsilon '} \leq
C 2^ { -\eta k}.
\end{equation}
Summing over $k$ will give the estimate for $ \nontan{(\nabla u )}$.

We begin
with the Cauchy-Schwarz inequality and then use the estimate of Lemma
\ref{63} (note that we have $\epsilon>2\beta/(\pi - 2\beta)>0>-1+2/q_0$)
 to obtain  
\begin{eqnarray*}
\int_{\Sigma_k } ( \nabla u ) ^ * (x) \, d\sigma _{\epsilon '} ( x) 
& \leq&
\left( \int _{ \Sigma _k }  \nontan{(\nabla u ) } (x) ^ 2 \, d\sigma _
\epsilon (x) \right) ^ { 1/2} 
\sigma_{ 2\epsilon ' -\epsilon } ( \Sigma _k ) ^{ 1/2} \\
 & \leq &
C  \frac { \sigma _{ 2\epsilon ' - \epsilon } ( \Sigma _k ) ^ { 1/2}
  \sigma _ \epsilon ( \Sigma _k ) ^ { 1/2} }{\sigma_{ \epsilon '} (
   \Delta _ \rho ( x _ a) )  2 ^ { k(1+  \delta )}} \\
 & = & C 2 ^ { -\delta k } \frac { \max ( |x_a| , 2^ k \rho ) ^ { \epsilon
    ' } }{ \max ( |x_a| , \rho ) ^ { \epsilon '} }\, 
%=:\,
%C2^ { -\delta k }M^{\epsilon'}_k, 
\end{eqnarray*}
where the last inequality follows from Lemma \ref{51}. 
Note that we have
\begin{equation}\label{Emax}
\frac { \max ( |x_a|, 2 ^ k \rho ) } { \max ( |x_a|, \rho ) } = \left\{
\begin{array}{ll} 1, \qquad & \frac{|x_a|}{\rho}<1 \\
|x_a|/\rho , \qquad& 1< \frac{|x_a|}{\rho} < 2^k \qquad \\
2^ k , \qquad  & \frac{|x_a|}{\rho}>2^ k .
\end{array}\right.
\end{equation}
In particular,  we have
$$ 1\leq\ \frac { \max ( |x_a|, 2 ^ k \rho ) } { \max ( |x_a|, \rho ) }
\ \leq\, 2^k\, .$$
Thus, letting $\eta := \delta - \mathrm{max}\{ \epsilon', 0\}>0$, 
we obtain  (\ref{64}). 
\comment{
\begin{equation}
%{eqnarray*}
\int_{\Sigma_k } ( \nabla u ) ^ * (x) \, d\sigma _{\epsilon '} ( x)
   \leq  C2 ^ { -\eta k }
%\left \{ \begin{array} {ll}
%2^ { \epsilon ' k } , \qquad & \epsilon ' > 0 \\
%1, \qquad   & \epsilon ' \leq 0 .
%      \end{array}
%\right.
%\end{eqnarray*}
\end{equation}
%by  Lemma \ref{51}. Note that we have 
%$$
%\frac { \max ( |y_a|, 2 ^ k \rho ) } { \max ( |y_a|, \rho ) } = \left\{
%\begin{array}{ll} 1, \qquad & |y_a| < \rho \\
%|y_a|/\rho , \qquad& \rho < |y_a| < 2^k \rho \\
%2^ k , \qquad  & 2^ k \rho < |y_a|.
%\end{array}\right. 
%$$

%Hence, since $ -\delta + \max (\epsilon ', 0) =-\delta < 0$, it follows that 
\noindent
% Thus, (\ref{64}) holds with $\eta = \delta $.
with $\eta>0$ (provided $\epsilon'<\delta$), as desired.
}
 Summing over $k$ yields estimate (\ref{EMPH11}).
% of Theorem \ref{53.5}.  

Now we indicate why $ \partial u / \partial \nu$ lies in $H^1 (
\partial\Omega, d\sigma_{ \epsilon '})$.
  A similar, but simpler, argument shows
that 
$u$ lies in $H^{1,1}(\partial\Omega, d\sigma _ {\epsilon '})$. 
We first prove the vanishing moment condition:
\begin{equation}\label{Emoment}
\int _{ \partial \Omega } \frac {\partial u }{ \partial \nu } \, d
\sigma = 0. 
\end{equation}
To this end, we 
 observe that, since $\int_{\partial\Omega}a\,d\sigma =0$,
 we may proceed
 as in the proof of Lemma \ref{62} and inequality (\ref{64})  
%except for the normalization,  a 
%$\sigma_{\epsilon '}$-atom is also a $ \sigma$-atom. Thus,
%the argument above  shows that 
to obtain
\begin{equation} \label{itsl1}
\int_{ \partial \Omega }  \nontan{(\nabla u )} \, d\sigma
<+\infty.  
\end{equation}
%The estimate (\ref{itsl1}) implies that 
%\begin{equation}\label{itsmean0}
%\int_{ \partial \Omega }  \frac {\partial u }{ \partial \nu } \,
%d\sigma = 0.
%\end{equation}
%To prove (\ref{itsmean0}),  we use
On account of (\ref{itsl1}) we may now apply the divergence theorem
 to obtain (\ref{Emoment}).
%and (\ref{itsl1})
%implies that  the boundary term at infinity vanishes. 

\noindent 
To continue, we follow the arguments of Coifman and Weiss \cite{CW:1976}.
% for 
%studying molecules to show that $ \partial u /\partial \nu$ lies in
%$H^1 ( \sigma_{\epsilon '})$.
  We begin by writing 
$$
\frac {\partial u }{ \partial \nu } = \sum_{ k = 0 }^ \infty b_k
$$
where $ b_0  :=  \chi_{ \Delta_0}\big(\frac {\partial u }{ \partial \nu }
 % -\frac { m_0}{ \sigma(\Delta _0)} \chi _{ \Delta _0}$ and
-\int_{\Delta _{\rho}}
\!\!\!\!\!\!\!\!\!\!-\ \ \frac {\partial u }{ \partial \nu }d\sigma\big)$
 and, for $ k\geq 1$,
$b_k :=  \chi_{ \Sigma_k}\frac {\partial u }{ \partial \nu }
+\chi _{ \Delta _{k-1}}
\int_{\Delta _{k-1}}
\!\!\!\!\!\!\!\!\!\!\!\!\!\!\!-\ \ \ \ 
\frac {\partial u }{ \partial \nu }d\sigma
-\chi _{ \Delta _{k}}
\int_{\Delta _{k}}
\!\!\!\!\!\!\!\!\!\!-\ \ \frac {\partial u }{ \partial \nu }d\sigma$.
%  +\frac { m_{k-1}}{ \sigma(\Delta _{k-1})} \chi _{ \Delta _{k-1}}
%  -\frac { m_k}{ \sigma(\Delta _k)} \chi _{ \Delta _k}$
% for $ k\geq 1$.
%In these definitions, we have put  $ m_k = \int _{ \Delta _k }
%\frac {\partial u }{   \partial \nu } \, d \sigma $. 
The estimate of Lemma  $\ref{63}$ and arguments used above imply that
$$
\int_{\partial \Omega} |b_k| \, d\sigma_{\epsilon'} \leq \left( \int_{ \Delta _k }
|b_k |^2 \, d\sigma _ \epsilon  \right )^ { 1/2} \sigma_{ 2\epsilon '
  - \epsilon } 
( \Delta _k ) ^{ 1/2}  \leq 2^ { -\eta k} .
$$
where again $ \eta := \delta - \max( \epsilon ' , 0 )>0$.   The one tricky
point in this argument is that we must use (\ref{Emoment}) to obtain
that $ \int _{ \Delta _k }
\frac {\partial u }{   \partial \nu } \, d \sigma =
 - \int _{ \partial \Omega \setminus \Delta _k } \frac {
\partial u }{ \partial \nu } \, d \sigma $.   Thus we have that $
2^{\eta k} b_k$ is normalized in
 $L^1( d \sigma _{\epsilon'})$.
% Recall
%that we 
%defined our atoms to lie in $L^ \infty$. Thus,
In order to prove that $b_k $ is in
$H^ 1 ( d\sigma _ {\epsilon '})$, one now proceeds
 as in Str{\"o}mberg and Torchinsky \cite[Theorem 1, p. 111]{ST:1989}
 to show that the grand maximal
function is in $L^1 (  d\sigma_{\epsilon '})$ 
and then find an atomic decomposition.
 This
argument gives estimate (\ref{EMPH12}) for $ \partial u / \partial \nu$
in $H^1 ( d \sigma_{ \epsilon '})$; the corresponding
 estimate for $u$ is obtained in a similar fashion. 
\end{proof}

\section{Existence and uniqueness for $H^1$ and $L^p$}
\setcounter{equation}{0}

In this section, we give the final arguments to prove our main result,
Theorem \ref{maintheorem} existence and uniqueness for the mixed
problem.
 
In Section 4 we obtained existence and uniqueness for the solution of the
$L^2(d\sigma_\epsilon)$-mixed problem for $ \epsilon$ in an interval
which includes positive values of $ \epsilon$ and 
does not include $0$. In this section, we will study solutions which
have atomic  data and show that the non-tangential maximal function
for such solutions will lie in $L^1( d\sigma_{ \epsilon'})$ for 
 $\epsilon'$ small. 
%$ p=1$. We will consider Hardy space results in $H^1 ( \sigma_{ 
%  \epsilon'})$ for $ \epsilon'$ small and negative.

% result for $H^1 (d\sigma_{ \epsilon'})$
%where $M$ denotes the Lipschitz constant
% of the domain. 
Our results are restricted to domains with Lipschitz constant $M$ 
 that is less than
1. This restriction is inherited from the previous
sections  (Lemma \ref{LVFMP}). 
%For such domains, we find an exponent $p_0=p_0(M)>1$ so that 
%the mixed problem is uniquely solvable in
% $L^p(\partial\Omega, d\sigma)$ for $ 1< p < p_0$.
 We do not make an effort to find the largest value of $p_0$ (nor do we
expect that the restriction $M<1$ is essential). 
 However, as the example discussed in the introduction indicates, 
 we cannot expect to always  have $p_0\ge 2$,  
 even in the case of smooth domains. 
%% \marginpar{Need ref?}
Furthermore, 
in Brown \cite{RB:1994b} the mixed problem is solved in
$L^2(d\sigma)$ for certain Lipschitz domains with arbitrarily large
Lipschitz constant; 
thus, it is  more  than the Lipschitz
constant that governs the solvability of this problem.

We begin with our uniqueness result. 

\begin{lemma}  \label{lpuniq} Let $ \delta >0 $ be as in  Lemma \ref{58.5}.
Suppose that  $\epsilon '$ and $p\geq 1 $ satisfy:
$-\delta <\epsilon ' \leq 0 $, $ 1/p'- \epsilon'/p  <
\delta$.
Under these hypotheses, the solution of the mixed problem
(\ref{Mixed}) for $L^p( d\sigma_{ \epsilon '})$ 
is unique. 
\end{lemma}

\begin{proof} Suppose that $u$ solves the mixed problem for $ L^p(
d\sigma_{ \epsilon '})$ with zero data, that is
\begin{equation}\label{EMP0}
\left\{
\begin{array}{ll}
\Delta u = 0 \qquad & \mbox{in } \Omega \\
u= 0 \qquad & \mbox{on } D \\
\frac { \partial u }{ \partial \nu } =
0 
& \mbox{on } N \\
(\nabla u )^* \in   L^p( d\sigma_{\epsilon'} )
\end{array}
\right.
\end{equation}
We will show that $u$ solves the regularity problem (\ref{Regular}) with zero
data and then use Lemma  \ref{RPUniq} to conclude that $u = 0$. To see
that $u$ vanishes on $N \subset \partial \Omega$, we will show that there is
$\eta>0$ such that for any $H^1(d\sigma_{\eta})$-atom, $a$, we have 
%argue by
%duality using the solutions of (\ref{Mixed}) which were constructed in
%Theorem \ref{53.5}. 
\begin{equation}\label{Eu0}
\int_N u\,a\,d\sigma =0\, .
\end{equation}
%for any $\sigma_{\eta}$-atom, $a$,
This means that is $u$ is constant almost everywhere.
In order to prove (\ref{Eu0}), we fix $ \eta > 0$ so that
  $$ 
0 < \frac{1}{p'}- \frac{\epsilon'}{p}  < \eta < \delta
 $$
(where $1/p+1/p'=1$) and we let $a$ be a $\sigma_{\eta}$-atom
supported in a boundary ball $\sball {x_a} r$. According to (the proof
of) Theorem \ref{53.5}, 
% there is a solution $v$ to
the mixed problem
% (\ref{Mixed})
in $H^ 1 ( d\sigma _{ \eta})$
 with Neumann data $a$ and with zero Dirichlet data has a solution
$v$ that satisfies
\begin{equation}
\label{V1}
|v(z) | \leq C  |z-x_a|^{ -\delta } , \quad
 \mbox{if\ \ $|z-x_a| \geq 2\rho$}\, ;
\end{equation}
 moreover, for $ \epsilon $  as in Theorem \ref{MixedProb2}, we have
\begin{equation}\label{V2}
\nontan{( \nabla v )} \in L^ 2 (  d\sigma _ \epsilon)\cap
 L^ 1 ( d\sigma _\eta).
\end{equation}
% for $ \epsilon $  as in Theorem \ref{MixedProb2}
% in Lemma \ref{58.5}. This
% implies that $v$ is bounded (see the argument below for $u$). The
% solution 
% $u$ satisfies  our hypothesis that $ \nontan { (\nabla u )}$ lies in
% $L^ p ( \sigma _{ \epsilon '})$. 
We will need a pointwise estimate for $u$; to this end, given $z$ in 
 $\Omega$, we consider the path in $\Omega$ from $0$ to $z$
 given by: $ \gamma_z ( t) = ( tz_1, \phi(tz_1)+ t ( z_ 2 -
 \phi(z_1)))$. Recall that $ \phi$ is the function whose graph gives $
 \partial \Omega$. We define the following Carleson measure with respect
 to $d\sigma$:
$$ d\mu_z (y) := \chi_{B_{2|z|}}(y)\, d\mathcal H_1\big|_{\gamma_z}\, ,
\quad
 y\in \Omega\, ,$$
where $\mathcal H_1$ denotes $1-$dimensional Hausdorff measure.
By the Fundamental Theorem of Calculus, properties of Carleson  measures
(Proposition \ref{Carleson}) and H\"older inequality, we have 
%Integrating along this path and using Proposition
% \ref{Carleson} on Carleson measures gives that 
\begin{eqnarray}
 |u(z)| & \leq&  \int_{\gamma_z}|\nabla u|\, |d\gamma_z|  \nonumber   \\ 
        & \leq&  C\int_{\Delta_{2|z|}(0)}(\nabla u)^*d\sigma
        \nonumber  \\
        & \leq &  |z|^{ \frac 1 {p'} - \frac { \epsilon ' }
          p } \| \nontan { ( \nabla u ) } \|_{ L^ p (d\sigma _{ \epsilon ' } 
  ) } .\label{growth}
\end{eqnarray}
Note that by a similar argument we may show that $u$ is locally bounded
 on $\partial\Omega$.

With these estimates collected, we now proceed to the main part of the
argument. 
 Let $R$ be large and let $ \psi_R$ be a cutoff function which is equal to 
1 on
$\ball 0 R $, zero outside $ \ball 0 {2R}$ and such that $|\nabla \psi_R|\leq
C/R$. We apply Green's second 
identity to the pointwise products $v\,\psi_R$ and $u\,\psi_R$ and obtain 
\begin{equation}\label{Green}
\int _{ \partial \Omega } \psi_R^2 \left( v \frac { \partial u }{ \partial \nu
} - u \frac { \partial v } { \partial \nu }\right)\, d\sigma  = 2
\int _ \Omega \psi_R\big(v \nabla u \cdot \nabla \psi _R -
 u \nabla v \cdot \nabla
\psi _R\big) \, dy .
\end{equation}
Concerning the left-hand side of (\ref{Green}), on account of the boundary
 conditions satisfied by $u$ and $v$, we have
$$
\int _{ \partial \Omega } \psi_R^2 \left( v \frac { \partial u }{ \partial \nu
} - u \frac { \partial v } { \partial \nu }\right)\, d\sigma  =
-\int_N\psi_R^2\, u\, a\, d\sigma\, .
$$
 Since $a$ is compactly supported and $u$ is locally bounded
 (see (\ref{growth}) and comment thereafter) we may apply the Lebesgue
 dominated convergence theorem and conclude
\begin{equation}\label{ELHS}
\lim\limits_{R\to \infty}
\int _{ \partial \Omega } \psi_R^2 \left( v \frac { \partial u }{ \partial \nu
} - u \frac { \partial v } { \partial \nu }\right)\, d\sigma  =
-\int_N u\, a\, d\sigma\, .
\end{equation}
%A routine limiting argument is needed to justify this. This argument
%uses that $\nontan{( \nabla u)}$ and  $\nontan{( \nabla v)}$ are in
%$L^ 1_{ loc}$  and that $u$ and $v$ are locally bounded.  
Concerning the right-hand side of (\ref{Green}),
we will show that
\begin{equation} \label{limit}
\lim_{ R\rightarrow \infty } \int _ \Omega \psi_R\,v \nabla u \cdot \nabla
\psi_R\, dy\ =\ 
\lim_{ R\rightarrow \infty } \int _ \Omega \psi_R\, u \nabla v \cdot
 \nabla  \psi_R  \, dy\   =\  0 .  
\end{equation}
Indeed, if $ R$ is large, then from Lemma \ref{58.5}, for $z$ in the
support of $ \nabla \psi_R$, we have
 \begin{equation}
\label{V3}
|v(z) | \leq C R^{ -\delta }\,  
\end{equation}
 and it follows  that 
$$
\int_ \Omega  \psi_R\, v \nabla u  \cdot \nabla \psi _R   \, dy 
\leq C R^{ -\delta } \int _\Omega | \nabla \psi _R | |\nabla u |\, dy
.
$$
By H\"older inequality we have
\begin{eqnarray*}
\int _\Omega | \nabla \psi _R | |\nabla u |\, dy &\leq &
R^{-\frac{\epsilon'}{p}}\!\!
\int_{B_R(0)}\! |\nabla\psi_R|\,|\nabla
u|\,|y|^{\frac{\epsilon'}{p}}dy \\
& \leq &
\comment{
R^{-\frac{\epsilon'}{p}}\left(\int_{B_R(0)}|\nabla u|^p|y|^{\epsilon'}dy\right)
^{\frac{1}{p}}\!\!\! R^{-1+\frac{2}{p'}}\,=\,
}
R^{\frac{1}{p'}-\frac{\epsilon'}{p}}\!\left(
\int_{\Omega}|\nabla u|^p\,\chi_{B_R(0)}\,\frac{|y|^{\epsilon'}}{R}dy\right)
^{\frac{1}{p}}.
\end{eqnarray*}
By Carleson Theorem, see (\ref{Ecarleson}) and Proposition \ref{Carleson},
the last term in the inequalities above is bounded by
$$
CR^{\frac{1}{p'}-\frac{\epsilon'}{p}}\!\|(\nabla u)^*\|
_{L^p( d\sigma_{\epsilon'})}\, .
%\left(\int_{\partial\Omega}
%((\nabla u)^*)^pd\sigma_{\epsilon'}\right)^{\frac{1}{p}}
$$
We conclude
\begin{equation}\label{Ezero1}
\left|\int_{\Omega}\psi_R(v\nabla u\cdot\nabla\psi_R)dy\right|\, \leq\,
CR^{\frac{1}{p'}-\frac{\epsilon'}{p}-\delta}
\|(\nabla u)^*\|
_{L^p( d\sigma_{\epsilon'})};
% \to 0\ \ \mathrm{as\ }\ R\to\infty\, ,
\end{equation}
%since $1/p'-\epsilon'/p-\delta <0$.
%$$
%\int _\Omega | \nabla \psi _R | |\nabla u |\, dx \leq C R ^ { 1/ p'-
%  \epsilon ' / p } 
%\left( \int_{\partial \Omega }{ \nontan{ ( \nabla u ) }} ^ p \, d\sigma
% _{ \epsilon '} \right ) ^ { 1/ p ' } .
%$$
our hypotheses on $ \eta$, $ \epsilon ' $ and $p$ imply that 
the exponent of $R$ is negative, so the first
 integral in (\ref{limit}) vanishes as $ R \rightarrow \infty$. 
To handle the second integral in (\ref{limit}) we apply 
(\ref{growth}); using the fact the $\grad \psi_R$ is supported in
the annulus $R<|z|<2R$ we obtain
$$
\left|\int_{\Omega}\psi_R\, u\nabla\psi_R\cdot\nabla v\, dy\right|\,\leq\,
C R ^  { \frac 1 {p'} - \frac { \epsilon ' } p }
\|(\nabla u)^*\|
_{L^p(\partial\Omega. d\sigma_{\epsilon'})}
R^{-\eta}\int_{R<|y|<2R}|\nabla v|\frac{|y|^{\eta}}{R}dy.
$$
By applying  Proposition \ref{Carleson} on Carleson measures one more
time, we see that the latter is bounded by
$$
%\int _\Omega u \nabla \psi _ R\cdot  \nabla v \, d x 
%\leq
  C R ^  { \frac 1 {p'} - \frac { \epsilon ' }{ p}-\eta  }
\|(\nabla u)^*\|
_{L^p( d\sigma_{\epsilon'})}
 \int _{\partial \Omega } \nontan { ( \nabla v ) } \,
d \sigma _ { \eta } .
$$
 It follows that the
second integral in (\ref{limit}) also vanishes as $ R$ tends to infinity. This
completes the proof of (\ref{limit}) and of this Lemma. 
\end{proof}

\begin{theorem} \label{TMPNA}  Suppose $\Omega$, $N$, $D$ is a standard graph
domain for the mixed problem, with Lipschitz constant $M<1$ and set $
\beta = \arctan(M)$. 
%Let $\delta>0$ be as in Lemma \ref{}. 
% satisfying> the conditions of Lemma
%\ref{LVFMP},
There is $\delta=\delta(M) $ satisfying $ 0< \delta <1$ so that, 
for
%for $2\beta/(\pi - 2\beta)<\epsilon<1$ 
%(where we have set $\beta = \tan^{-1}M>0$) and
$\max \{-\delta\, , \frac{4\beta-\pi}{2(\pi-\beta)}\}<\epsilon'\le 0 $
%(where we have set $\beta = \tan^{-1}M>0$),
% there is an $\epsilon_0 = \epsilon_0 (M)$ so that, if $ |\epsilon | <
%\epsilon_0$ and $a$ is an $H^1 (N, d\sigma_{ \epsilon'})$-atom,
%we may uniquely solve
the mixed problem  (\ref{Mixed}) for
$L^1( d\sigma_{\epsilon'})$ with
Dirichlet data in $H^{1,1}(D, d\sigma_{\epsilon'})$ and
with Neumann data in
$H^1 (N, d\sigma_{ \epsilon'})$
% $a$
% and with zero
%Dirichlet data.
is uniquely solvable.
 The solution $u$ satisfies the following
estimates
\begin{equation}\label{EMPhardy}
\int_{ \partial \Omega} (\nabla u )^* (x) \, d \sigma_{ \epsilon'}
\leq C(M, \epsilon')\, \left(
\|f_N \| _{ H^{1}(N, d\sigma_{\epsilon'})} +
   \| f_D \|_
{ H^{1,1} (D, d\sigma_\epsilon')}\right) ;
\end{equation}
\begin{equation}\label{EMPnormal}
\!\!\!\!\!\!\!\!
\|u \| _{ H^{1,1}( d\sigma_{\epsilon'})} +
   \left\|\frac {\partial u } {\partial \nu}\right\|_{ H^1 (d\sigma_{\epsilon'})}\!\!\!\!\!\!\!\!\!\!\!\!\!\!  \le\,
   C(M, \epsilon')\left(
\|f_N \| _{ H^{1}(N, d\sigma_{\epsilon'})} +
 \| f_D \|_
{ H^{1,1} (D, d\sigma_\epsilon')}\right).
\end{equation}
\end{theorem}

\begin{proof}
We first observe that
we may use the solution of the regularity
problem from Theorem \ref{rnh1} to reduce to the case where the
Dirichlet data is zero. More precisely,
% if
% $f_D \in H^{1,1}(D, d\sigma_{\epsilon'})$ denotes the Dirichlet data for
% the mixed problem,
 we consider the (unique) solution $v$ of the regularity
 problem (\ref{Regular}) with data $\tilde f_D \in H^{1,1}( d\sigma_{\epsilon'})$
 (here $\tilde f_D$ denotes an extension of $f_D$ to $\partial\Omega$).
By Theorem \ref{rnh1} it follows that
 $\partial v/\partial\nu \in H^1(\partial\Omega, d\sigma_{\epsilon'})$,
 and it is easy to see that $u$ is the unique solution of
 the mixed problem with data
$f_D$ and $f_N$ if and only if $u-v$ solves the mixed problem
with zero Dirichlet data and with Neumann data
 $g_N:=f_N - \partial v/\partial\nu$.

 To prove existence 
 we consider any atomic decomposition for $g_N$, namely
\begin{equation}\label{atomicgN}
g_N(x) = \sum _{ j=1}^ \infty \lambda _j a _ j (x),\quad 
\sum _{ j=1}^ \infty |\lambda _j|<+\infty\, ,
\end{equation}
see (\ref{EH1}) and (\ref{EH11}). By Theorem \ref{53.5} it follows that
for each $j$ the mixed problem:
\begin{equation}\label{mixedj}
\left\{
\begin{array}{ll}
\Delta h_j = 0 \qquad & \mbox{in } \Omega \\
h_j= 0 \qquad & \mbox{on } D \\
\frac { \partial h_j }{ \partial \nu } =
 a_j
& \mbox{on } N
 \\
(\nabla h_j )^* \in   L^1( d\sigma_{\epsilon'} )
\end{array}
\right.
\end{equation}
has a solution $h_j$ that satisfies:
\begin{equation}\label{Mixedjest}
 \|h_j \| _{ H^{1,1}( d\sigma_{\epsilon'})} +
   \left\|\frac {\partial h_j }
 {\partial \nu}\right\|_{ H^1 ( d\sigma_{\epsilon '} )}
   \leq C(M, \epsilon').
\end{equation} 
Thus, the function
$$
h:= \sum _{ j=1}^ \infty \lambda _j h_ j 
$$
is a solution of the mixed problem
with zero Dirichlet data and with Neumann data
 $g_N$, and it satisfies:
$$
\left\|\frac{\partial h}{\partial \nu} \right\|_
{H^1( d\sigma_{\epsilon'})} \, \leq\,
C(M,\epsilon')\sum_{j=1}^{\infty}|\lambda_j|\, . 
$$
By Lemma \ref{lpuniq}, $h$ is unique, (i.e. $h$ is independent of
 the choice of the atomic decomposition for $g_N$); taking the infimum over
 all atomic decompositions of $g_N$ now yields (\ref{EMPhardy}).
This proves existence; uniqueness follows from Lemma \ref{lpuniq}.  
\end{proof}
Next, we recall a few well known results concerning the
the complex  interpolation spaces  of weighted $L^p$ and Hardy spaces, see
Bergh and L\" ofstrom \cite[Theorem 5.5.3, Corollary 5.5.4]{BL:1976},
Str\"omberg and Torchinsky \cite[Theorem 3, pg.~179]{ST:1989}.
%, Taylor \cite[(Proposition 2.1)]{MT:1996}.
%%Does this reference discuss weighted spaces? 
For weights $ w_0$ and $w_1$, and  exponents $p_0, p_1$ we
set 
\begin{equation}\label{Einterp1}
\frac 1 {p_{\theta}} =  \frac { 1-\theta } {p_0} + \frac \theta {p_1},\quad
w_\theta = w_0^{\frac {(1-\theta)p_\theta}{p_0}}w_1^{\frac
  {\theta p_\theta}{p_1} }, \quad 0\leq \theta \leq 1. 
\end{equation}
We let $ [A,B]_\theta$ denote the complex interpolation space of index
$ \theta$ as defined in the monograph of Bergh and L\"ofstrom
\cite[Chapt.~4]{BL:1976}.  
\begin{theorem} \label{70}
Suppose $ w_{0} $ and $ w_{ 1}$ are weights, then we have
$$
 [L^1 ( w_0\, d\sigma),
L^2 (  w_1 \, d\sigma)]_\theta
= L^{p_\theta} (  w_\theta \, d\sigma).
$$
If in addition the weights $w_j$ are in  $A_\infty( d\sigma)$ for $j=0,1$,
then we have 
$$
[ H^1 ( w_0\, d\sigma), 
L^2 (  w_1 \, d\sigma)]_\theta = L^ {p_\theta} (  w_\theta\,  d\sigma). 
$$

% Then, with same notations as (\ref{Einterp1})
% and \cite[Proposition 2.1]{MT:1996}, 
% we have
%\begin{eqnarray*}\label{Einterp2}  
%\end{eqnarray*} 
Moreover, if a linear operator $T$ is bounded:
\begin{eqnarray*}\label{Einterp3} 
T: H^1 ( w_0\, d\sigma)
&\to& L^1 ( w_0\, d\sigma), \\ 
T: L^ 2 (  w_1 \, d\sigma)
 &\to& L^ 2 (  w_1 \, d\sigma)
\end{eqnarray*}
 with norms $M_0$ and $M_1$ respectively, then
 $T$ is bounded:
$$L^ {p_\theta} (  w_\theta\, d\sigma)\to
 L^ {p_\theta} ( w_\theta \,  d\sigma)$$
with norm $M$ satisfying
$$
M\leq C M_0^{1-\theta}M_1^{\theta}. 
$$
\end{theorem}

\begin{remark} The constant $C$  in the estimate for the operator norm is
1 when we consider Lebesgue spaces. It may not be one for Hardy
  spaces, see Str\"omberg and Torchinsky \cite{ST:1989}. 
\end{remark}

%%More details?? RMB 8-16-04
We will focus on the case:
% $p_0=1$; 
$w_0 \, d\sigma= d\sigma_{\epsilon'}$, with $\epsilon'<0$ as in
 Theorem \ref{TMPNA}, and 
%$p_1=2$;
 $w_1\,d\sigma= 
d\sigma_{\epsilon}$, where
 $\epsilon>0$ is as in Theorem \ref{MixedProb2}.
%Since odd reflection extends $H^1(N; d\sigma_{\epsilon '}) $ to $H^1 (
%\partial \Omega ; d\sigma_{ \epsilon '} )$, it is easy to generalize the
%above result to spaces defined on $N$.
We are now ready to give the proof of our  main result, Theorem
\ref{maintheorem}. 

\comment{
\begin{theorem} Let $ \Omega $, $N$, and $D$, be a standard graph
domain for the mixed problem (\ref{Mixed})
 with Lipschitz constant $M < 1$. 
% and let
%$0<\delta =\delta(M)<1$ be as in (\ref{54.3}).
% (so that
%Theorem \ref{MixedProb2} applies).
 Then, there is a value $ p_0>1$, $p_0=p_0(M)$ so that
for $ 1 < p < p_0$, the (unweighted) mixed problem for
 $L^ p (d\sigma)$ is uniquely solvable.
 The solution $u$ satisfies
$$
\| ( \nabla u )^ *\|_{ L^ p ( d \sigma)}\leq C(p,M)
 \left( \| f_N
  \|_{L^p(N, d\sigma)} + \left\|\frac{\partial f_D}{d\sigma}\right\| _{L^ p (D, d\sigma)}\right).
$$
\end{theorem}
}
\begin{proof}[Proof of Theorem \ref{maintheorem}]
% To establish existence, we apply Theorems
%\ref{MixedProb2} and \ref{TMPNA} 
% and for $0<\delta<1$ as in
%(\ref{54.3}), we fix 
We first use  a  result of Dahlberg and Kenig \cite[Theorem 3.8]{DK:1987} to
reduce to the case  where the Dirichlet data in the mixed problem is zero. 
(While Dahlberg and Kenig only discuss $n\geq 3$ in their work, one
can  extend their results to two dimensions.)  
%%X12A 

%To establish existence, we 
We will use interpolation to  establish existence of solutions
satisfying the estimate (\ref{EMPLP}). Because the
complex method applies to linear operators, we employ a standard
technique to obtain the non-tangential maximal function as a supremum
of linear operators.  Fix
$\{y_j\}_{j\in\naturals }$, a dense subset of the sector  $\Gamma(0)$
 and let 
$\mathcal E = \{E_j\}$ be decomposition of  $\partial\Omega$ into
disjoint measurable subsets,  $E_j$. We define a linear operator:
$$ f_N \in H^1(d\sigma_{\e'})+L^2(d\sigma_{\e})\ 
\to \  T_{\mathcal E}(f_N)(x) \assign  
\sum_j\chi_{E_j}(x)\nabla u (x+y_j)\, ,\quad x\in \partial\Omega,
$$
where $u$ is the solution of the $H^1(d\sigma_{\epsilon'})$-mixed problem
  (resp.~the
$L^2(d\sigma_{\e})$-mixed problem) for data  $f_N$ and $f_D=0$.
 Note that for a suitable sequence
of decompositions $\{\mathcal E_k\}$, we have
$$(\nabla u)^*(x) = \lim\limits_{k\to \infty} |T_{\mathcal E_k}f_N(x)|\quad
 a.e.\ \, x \in \partial \Omega\, .$$

Now complex interpolation, (see Theorem \ref{70}) implies that the
operator $T_{\cal E} $ is bounded on the intermediate spaces with a
norm independent of $ {\cal E}$. Fatou's lemma then yields 
 boundedness for the non-tangential maximal function. 
It is easy to see that for the spaces $L^2 ( d\sigma_ \epsilon )$, $
(2\beta)/(\pi - 2\beta)< \epsilon < 1$ and  $H^1 ( \sigma_{\epsilon '})$, $
\max(-\delta, ( \epsilon -1)/2 ) < \epsilon '\leq 0$, the intermediate
spaces include $ L^p( d\sigma)$ for $ 1< p < p_1(M)$ where
$$
p_1(M) = 
\frac 
{ 2 ( \pi-2\beta) \min( \delta, 
  \frac {\pi - 4 \beta}
	{2 ( \pi -\beta)}) + 2\beta }
{(\pi -2 \beta ) 
\min ( \delta , \frac { \pi - 4 \beta}{2(\pi -\beta)} )+2\beta
}.  
$$
In addition, in Lemma \ref{lpuniq}, uniqueness was established 
for $p$ in the range:
% (for $\epsilon'=0$):
 $1<p<p_2(M)$, $p_2(M)=1/(1-\delta(M))$. Thus, we have existence and
 uniqueness for $p \in ( 1, p_0 (M))$ with $ p_0 = \min(p_1(M),
 p_2(M))$. 
\comment{
To establish existence, we now apply  Theorem \ref{70} to 
interpolate between the $L^ 2 ( d\sigma_ \epsilon )$-result
 in Theorem \ref{MixedProb2} 
($2 \beta/( \pi -2 \beta) <\epsilon <1$) and the $H^ 1
(N,d\sigma_{\epsilon '})$-result in
Theorem \ref{TMPNA} (max$\{-\delta\, , (\epsilon -1)/2\}<\epsilon'< 0, $),  
%To find the upper bound $p_0$, 
where we choose $ \theta$ so that
$(1-\theta)\epsilon' + \theta \epsilon =0$ (that is, $w_{\theta} =1$).
  The corresponding value of $p$ (see (\ref{Einterp1})) is given by
$ p = 2( \epsilon -\epsilon')/ ( 2 \epsilon - \epsilon')$. 
For $\epsilon$ and $\epsilon'$ as above, the upper bound for $p$
occurs at 
$$p_1(M) = \frac{4\beta+
2\min\,\{\delta,\frac{\pi -4\beta}{2(\pi -\beta)}\}(\pi-2\beta)}
{4\beta+
\min\,\{\delta,\frac{\pi -4\beta}{2(\pi -\beta)}\}(\pi-2\beta)}.$$ 
Uniqueness for the unweighted mixed problem in $L^p( d\sigma)$
was established in Lemma \ref{lpuniq} for $p$ in the range:
% (for $\epsilon'=0$):
 $1<p<p_2(M)$, $p_2(M)=1/(1-\delta(M))$.
Thus, existence and uniqueness hold for $1<p<p_0=p_0(M)$ with
 $p_0 :=$ min$\{p_1, p_2\}$. 
}
\end{proof} 

%\marginpar{No it is not!}

%\marginpar{Oops--need to do uniqueness in $L^p$ and existence for $H^1$.) }

%\theequation

\section{Final remarks}
\setcounter{equation}{0}
We close by listing a few open questions.

\begin{itemize} \item  Can we remove the restriction that the
  Lipschitz constant of the domain is at most one?

\item Can we obtain similar results in higher dimensions?  The obvious
problem here is that our two-dimensional weighted estimates rely on
a Rellich identity which is based on complex function theory. 

\item Can we study domains where the boundary between $D$ and $N$ is
more interesting? For example in $\reals ^3$, let $ \Omega = \{ x: x_3
> c (|x_1|+|x_2|)\}$ and let $ D = \partial \Omega \cap \{ x: x_1x_2>
0\}$
and then $N= \partial \Omega \setminus D$. Can we solve the mixed
problem for some $L^p$ space in this domain? 

\item Can we extend the solution of the mixed problem to general
domains, rather than only graph domains? 
\end{itemize}

%\bibliographystyle{plain}
%\bibliography{main,inverse} 

\small \noindent \today 
\end{document}